\documentclass[reqno,10pt, centertags]{amsart}
\usepackage{amsmath,amsthm,amscd,amssymb,latexsym,upref}
\usepackage{esint,color}
\usepackage{hyperref}

\newcommand*{\mailto}[1]{\href{mailto:#1}{\nolinkurl{#1}}}
\newcommand{\arxiv}[1]{\href{http://arxiv.org/abs/#1}{arXiv:#1}}
\makeatletter
\def\theequation{\@arabic\c@equation}

\newcommand{\bbN}{{\mathbb{N}}}
\newcommand{\bbR}{{\mathbb{R}}}

\newcommand{\bbC}{{\mathbb{C}}}

\newcommand{\cB}{{\mathcal B}}

\newcommand{\cD}{{\mathcal D}}
\newcommand{\cE}{{\mathcal E}}

\newcommand{\cH}{{\mathcal H}}

\newcommand{\cK}{{\mathcal K}}

\newcommand{\cR}{{\mathcal R}}
\newcommand{\cS}{{\mathcal S}}

\newcommand{\cX}{{\mathcal X}}

\newcommand{\abs}[1]{\left\lvert#1\right\rvert}
\newcommand{\norm}[1]{\left\Vert#1\right\Vert}
\newcommand{\no}{\nonumber}
\newcommand{\lb}{\label}
\newcommand{\f}{\frac}

\newcommand{\ol}{\overline}
\newcommand{\ti}{\tilde}
\newcommand{\wti}{\widetilde}
\newcommand{\eps}{\varepsilon}

\newcommand{\ga}{\gamma}
\newcommand{\la}{\lambda}
\newcommand{\si}{\sigma}
\newcommand{\Om}{\Omega}
\newcommand{\dOm}{{\partial\Omega}}

\newcommand{\loc}{\text{\rm{loc}}}
\newcommand{\tr}{\text{\rm{tr}}}

\newcommand{\ran}{\text{\rm{ran}}}

\newcommand{\dom}{\text{\rm{dom}}}

\newcommand{\supp}{\text{\rm{supp}}}

\newcommand{\bi}{\bibitem}
\newcommand{\sgn}{\text{\rm{sign}}}
\newcommand{\hatt}{\widehat}
\renewcommand{\Re}{\text{\rm Re}}
\renewcommand{\Im}{\text{\rm Im}}
\renewcommand{\ln}{\text{\rm ln}}

\numberwithin{equation}{section}

\newtheorem{theorem}{Theorem}[section]
\newtheorem{lemma}[theorem]{Lemma}
\newtheorem{corollary}[theorem]{Corollary}
\theoremstyle{definition}
\newtheorem{definition}[theorem]{Definition}
\newtheorem{hypothesis}[theorem]{Hypothesis}
\newtheorem{remark}[theorem]{Remark}

\begin{document}

\title[Non-Self-Adjoint Operators and Infinite
Determinants]{Non-self-adjoint operators, Infinite Determinants, and
some Applications}
\author[F.\ Gesztesy, Y.\ Latushkin, M.\ Mitrea, and M.\
Zinchenko]{Fritz Gesztesy, Yuri Latushkin, \hspace{180pt} Marius
Mitrea, and Maxim Zinchenko \hspace{80pt}}

\address{Department of Mathematics,
University of Missouri, Columbia, MO 65211, USA}
\email{fritz@math.missouri.edu}
\urladdr{http://www.math.missouri.edu/personnel/faculty/gesztesyf.html}
\address{Department of Mathematics, University of
Missouri, Columbia, MO 65211, USA}
\email{yuri@math.missouri.edu}
\urladdr{http://www.math.missouri.edu/personnel/faculty/latushkiny.html}
\address{Department of Mathematics, University of
Missouri, Columbia, MO 65211, USA}
\email{marius@math.missouri.edu}
\urladdr{http://www.math.missouri.edu/personnel/faculty/mitream.html}
\address{Department of Mathematics, University of
Missouri, Columbia, MO 65211, USA}
\email{maxim@math.missouri.edu}
%%%%%%%%%%%%%%%%%%%%%%%%%%%%%%%%%%%%%%%%%%   
\dedicatory{With great respect and deep admiration, we dedicate this
paper\\ to the memory of Boris M.\ Levitan 1914--2004.}
%%%%%%%%%%%%%%%%%%%%%%%%%%%%%%%%%%%%%%%%% 
%\dedicatory{}
%\date{\today}
%\date{, 2003.}
\thanks{Based upon work supported by the US National Science
Foundation under Grant Nos.\ DMS-0405526, DMS-0338743, DMS-0354339,
DMS-0400639, FRG-0456306, and the CRDF grant
UP1-2567-OD-03.}
\thanks{\it Russ. J. Math. Phys. {\bf 12}, 443--471 (2005).}

%%%%%%%%%%%%%%%%%%%%%%%%%%%%%%%%%%%%%%%%% 
\begin{abstract}
We study various spectral theoretic aspects of non-self-adjoint operators.
Specifically, we consider a class of factorable non-self-adjoint
perturbations of a given unperturbed non-self-adjoint operator and provide
an in-depth study of a variant of the  Birman--Schwinger principle as well
as local and global Weinstein--Aronszajn formulas.

Our applications include a study of suitably symmetrized (modified)
perturbation determinants of Schr\"odinger operators in dimensions
$n=1,2,3$ and their connection with Krein's spectral shift function
in two- and three-dimensional scattering theory. Moreover, we study
an appropriate multi-dimensional analog of the celebrated formula by Jost
and Pais that identifies Jost functions with suitable Fredholm
(perturbation) determinants and hence reduces the latter to simple
Wronski determinants.
\end{abstract}
%%%%%%%%%%%%%%%%%%%%%%%%%%%%%%%%%%%%%%%%% 

\maketitle

{\scriptsize{\tableofcontents}}
%\normalsize

%%%%%%%%%%%%%%%%%%%%%%%%%%%%%%%%%%%%%%%%%
%%%%%%%%%%%%%%%%%%%%%%%%%%%%%%%%%%%%%%%%%
\section{Introduction} \lb{s1}
%%%%%%%%%%%%%%%%%%%%%%%%%%%%%%%%%%%%%%%%%
%%%%%%%%%%%%%%%%%%%%%%%%%%%%%%%%%%%%%%%%%

This paper has been written in response to the increased demand of
spectral theoretic aspects of non-self-adjoint operators in contemporary
applied and mathematical physics. What we have in mind, in particular,
concerns the following typical two scenarios: First, the construction of
certain classes of solutions of a number of completely integrable
hierarchies of evolution equations by means of the inverse scattering
method, for instance, in the context of the focusing nonlinear
Schr\"odinger equation in $(1+1)$-dimensions, naturally leads to
non-self-adjoint Lax operators. Specifically, in the particular case of
the focusing nonlinear Schr\"odinger equation the corresponding Lax
operator is a  non-self-adjoint one-dimensional Dirac-type operator.
Second, linearizations of nonlinear partial
differential equations around steady state and solitary-type solutions,
frequently, lead to a linear non-self-adjoint spectral problem. In the
latter context, the use of the so-called Evans function (an analog of the
one-dimensional Jost function for Schr\"odinger operators) in the course
of a linear stability  analysis has become a cornerstone of this circle of
ideas. As shown in \cite{GLM06}, the Evans function equals a (modified)
Fredholm determinant associated with an underlying Birman--Schwinger-type
operator. This observation naturally leads to the second main theme
of this paper and a concrete application to non-self-adjoint operators,
viz., a study of properly symmetrized (modified) perturbation
determinants of non-self-adjoint Schr\"odinger operators in dimensions
$n=1,2,3$.

Next, we briefly summarize the content of each section. In Section
\ref{s2}, following the seminal work of Kato \cite{Ka66} (see also Konno
and Kuroda \cite{KK66} and Howland \cite{Ho70}), we consider a  class of
factorable non-self-adjoint perturbations, formally given by
$B^* A$, of a given unperturbed non-self-adjoint operator $H_0$ in a
Hilbert space $\cH$ and introduce a densely defined, closed linear
operator $H$ in $\cH$ which represents an extension of $H_0+B^*A$. Closely
following Konno and Kuroda \cite{KK66}, we subsequently derive a general
Birman--Schwinger principle for $H$ in Section \ref{s3}. A variant of the
essential spectrum of $H$ and a local Weinstein--Aronszajn formula is
discussed in Section \ref{s4}. The corresponding global
Weinstein--Aronszajn formula in terms of modified Fredholm determinants
associated with the Birman--Schwinger kernel of $H$ is the content of
Section \ref{s5}. Both, Sections \ref{s4} and \ref{s5} are modeled after
an exemplary treatment of these topics by Howland \cite{Ho70} in
the case where $H_0$ and $H$ are self-adjoint. In Section \ref{s6} we
turn to concrete applications to properly symmetrized (modified)
perturbation determinants of non-self-adjoint Dirichlet- and
Neumann-type Schr\"odinger operators in $L^2(\Om;d^n x)$ with
$\Om=(0,\infty)$ in the case $n=1$ and rather general open domains
$\Om\subset \bbR^n$ with a compact boundary in dimensions $n=2,3$. The
corresponding potentials $V$ considered are of the form $V\in
L^1((0,\infty);dx)$ for
$n=1$ and $V\in L^2(\Om;d^nx)$ for $n=2,3$. Our principal result in
this section concerns a reduction of the Fredholm determinant of the
Birman--Schwinger kernel of $H$ in $L^2(\Om;d^n x)$ to a Fredholm
determinant associated with operators in $L^2(\partial\Om;
d^{n-1}\sigma)$. The latter should be viewed as a proper
multi-dimensional extension of the celebrated result by Jost and Pais
\cite{JP51} concerning the equality of the Jost function (a Wronski
determinant) and the associated Fredholm determinant of the underlying
Birman--Schwinger kernel. In Section \ref{s7} we briefly discuss an
application to scattering theory in dimensions $n=2,3$ and re-derive a
formula for the Krein spectral shift function (related to the logarithm
of the determinant of the scattering matrix) in terms of modified
Fredholm determinants of the underlying Birman--Schwinger kernel. We
present an alternative derivation of this formula originally due to
Cheney \cite{Ch84} for $n=2$ and Newton \cite{Ne77} for $n=3$ (in the
latter case we obtain the result under weaker assumptions on the
potential $V$ than in \cite{Ne77}). Finally, Appendix \ref{sA}
summarizes results on Dirichlet and Neumann Laplacians in
$L^2(\Om;d^nx)$ for a general class of open domains
$\Omega\subseteq\bbR^n$, $n\geq 2$, with a compact boundary. We prove the
equality of two natural definitions of Dirichlet and Neumann Laplacians
for such domains and prove mapping properties between appropriate scales
of Sobolev spaces. These results are crucial ingredients in our treatment
of modified Fredholm determinants in Section \ref{s6}, but they also
appear to be of independent interest.

We will use the following notation in this paper. Let $\cH$ and $\cK$
be separable complex Hilbert spaces, $(\cdot,\cdot)_{\cH}$ and
$(\cdot,\cdot)_{\cK}$ the scalar products in $\cH$ and $\cK$ (linear in
the second factor), and
$I_{\cH}$ and $I_{\cK}$ the identity operators in $\cH$ and $\cK$,
respectively. Next, let $T$ be a closed linear operator from
$\dom(T)\subseteq\cH$ to $\ran(T)\subseteq\cK$, with $\dom(T)$
and $\ran(T)$ denoting the domain and range of $T$. The closure of a
closable operator $S$ is denoted by $\ol S$. The kernel (null space) of $T$
is denoted by
$\ker(T)$. The spectrum and resolvent set of a closed linear operator in
$\cH$ will be denoted by
$\sigma(\cdot)$ and $\rho(\cdot)$. The Banach spaces of bounded
and compact linear operators in $\cH$ are denoted by $\cB(\cH)$ and
$\cB_\infty(\cH)$, respectively. Similarly, the Schatten--von Neumann
(trace) ideals will subsequently be denoted by $\cB_p(\cH)$,
$p\in\bbN$. Analogous notation $\cB(\cH_1,\cH_2)$,
$\cB_\infty (\cH_1,\cH_2)$, etc., will be used for bounded, compact, etc.,
operators between two Hilbert spaces $\cH_1$ and $\cH_2$. In addition,
$\tr(T)$ denotes the trace of a trace class operator $T\in\cB_1(\cH)$ and
$\det_{p}(I_{\cH}+S)$ represents the (modified) Fredholm determinant
associated with an operator $S\in\cB_p(\cH)$, $p\in\bbN$ (for $p=1$ we
omit the subscript $1$). Moreover,
$\cX_1 \hookrightarrow \cX_2$ denotes the continuous imbedding of the
Banach space $\cX_1$ into the Banach space $\cX_2$.

Finally, in Sections \ref{s6} and \ref{s7} we will introduce various
operators of multiplication, $M_f$, in $L^2(\Om;d^nx)$ by elements $f \in
L^1_{\loc}(\Om;d^nx)$, where $\Omega\subseteq\bbR^n$ is open and nonempty.

%%%%%%%%%%%%%%%%%%%%%%%%%%%%%%%%%%%%%%%%
%%%%%%%%%%%%%%%%%%%%%%%%%%%%%%%%%%%%%%%% 
\section{Abstract Perturbation Theory} \lb{s2}
%%%%%%%%%%%%%%%%%%%%%%%%%%%%%%%%%%%%%%%%
%%%%%%%%%%%%%%%%%%%%%%%%%%%%%%%%%%%%%%%%

In this section, following Kato \cite{Ka66}, Konno and Kuroda
\cite{KK66}, and Howland \cite{Ho70}, we consider a class
of factorable non-self-adjoint perturbations of a given unperturbed
non-self-adjoint operator. For reasons of completeness we will
present proofs of many of the subsequent results even though most of
them are only slight deviations from the original proofs in the
self-adjoint context.

We start with our first set of hypotheses.

%%%%%%%%%%%%%%%%%%%%%%%%%%%%%%%%%%%%%%%
\begin{hypothesis} \lb{h2.1}
$(i)$ Suppose that $H_0\colon\dom(H_0)\to\cH$,
$\dom(H_0)\subseteq\cH$ is a densely defined, closed, linear operator
in $\cH$ with nonempty resolvent set,
\begin{equation}
\rho(H_0)\neq\emptyset, \lb{2.1}
\end{equation}
$A\colon\dom(A)\to\cK$, $\dom(A)\subseteq\cH$ a densely defined,
closed, linear operator from $\cH$ to $\cK$, and
$B\colon\dom(B)\to\cK$, $\dom(B)\subseteq\cH$ a densely
defined, closed, linear operator from $\cH$ to $\cK$ such that
\begin{equation}
\dom(A)\supseteq\dom(H_0), \quad \dom(B)\supseteq\dom(H_0^*).
\lb{2.2}
\end{equation}
In the following we denote
\begin{equation}
R_0(z)=(H_0-zI_{\cH})^{-1}, \quad z\in \rho(H_0). \lb{2.3}
\end{equation}
$(ii)$ For some (and hence for all) $z\in\rho(H_0)$, the operator
$-AR_0(z)B^*$, defined on $\dom(B^*)$, has a
bounded extension in $\cK$, denoted by $K(z)$,
\begin{equation}
K(z)=-\ol{AR_0(z)B^*} \in\cB(\cK). \lb{2.4}
\end{equation}
$(iii)$ $1\in\rho(K(z_0))$ for some $z_0\in \rho(H_0)$.
\end{hypothesis}
%%%%%%%%%%%%%%%%%%%%%%%%%%%%%%%%%%%%

That $K(z_0)\in\cB(\cK)$ for some $z_0\in\rho(H_0)$ implies
$K(z)\in\cB(\cK)$ for all $z\in\rho(H_0)$ (as mentioned in
Hypothesis \ref{h2.1}\,(ii)) is an immediate consequence of
\eqref{2.2} and the resolvent equation for $H_0$.

We emphasize that in the case where $H_0$ is self-adjoint, the
following results in Lemma \ref{l2.1}, Theorem \ref{t2.3}, and
Remark \ref{r2.4} are due to Kato \cite{Ka66} (see also \cite{Ho70},
\cite{KK66}). The more general case we consider here requires only
minor modifications. But for the convenience of the reader we will
sketch most of the proofs.

%%%%%%%%%%%%%%%%%%%%%%%%%%%%%%%%%%%%% 
\begin{lemma} \lb{l2.1}
Let $z,z_1,z_2\in\rho(H_0)$. Then  Hypothesis \ref{h2.1} implies the
following facts:
\begin{align}
& AR_0(z)\in\cB(\cH,\cK), \quad \ol{R_0(z)B^*}=[B(H_0^*-\ol z)^{-1}]^*\in
\cB(\cK,\cH), \lb{2.5}
\\ & \ol{R_0(z_1)B^*}-\ol{R_0(z_2)B^*}
=(z_1-z_2)R_0(z_1)\ol{R_0(z_2)B^*} \lb{2.6} \\
& \hspace*{3.35cm} =(z_1-z_2)R_0(z_2)\ol{R_0(z_1)B^*}, \lb{2.7} \\
& K(z)=-A\ol{[R_0(z)B^*]}, \quad K(\ol z)^*=-B\ol{[R_0(\ol z)^*A^*]},
\lb{2.8} \\
& \ran(\ol{R_0(z)B^*})\subseteq \dom(A), \quad
\ran(\ol{R_0(\ol z)^*A^*})\subseteq \dom(B), \lb{2.9} \\
& K(z_1)-K(z_2)=(z_2-z_1)AR_0(z_1)\ol{R_0(z_2)B^*} \lb{2.10} \\
& \hspace*{2.25cm} =(z_2-z_1)AR_0(z_2)\ol{R_0(z_1)B^*}. \lb{2.11}
\end{align}
\end{lemma}
%%%%%%%%%%%%%%%%%%%%%%%%%%%%%%%%%%%%% 
\begin{proof}
Equations \eqref{2.5} follow from the relations in
\eqref{2.2} and the Closed Graph Theorem. \eqref{2.6} and \eqref{2.7}
follow from combining
\eqref{2.5} and the resolvent equation for $H_0^*$. Next, let
$f\in\dom(B^*)$, $g\in\dom(A^*)$, then
\begin{equation}
(\ol{R_0(z)B^*}f,A^*g)_{\cH}=(R_0(z)B^*f,A^*g)_{\cH}
=(AR_0(z)B^*f,g)_{\cK}=-(K(z)f,g)_{\cK}. \lb{2.12}
\end{equation}
By continuity this extends to all $f\in\cK$. Thus,
$-A\ol{[R_0(z)B^*]}f$ exists and equals $K(z)f$ for all $f\in\cK$.
This proves the first assertions in \eqref{2.8} and \eqref{2.9}. The
remaining assertions in \eqref{2.8} and \eqref{2.9} are of course
proved analogously. Multiplying \eqref{2.6} and \eqref{2.7} by $A$
from the left and taking into account the first relation in
\eqref{2.8}, then proves
\eqref{2.10} and \eqref{2.11}.
\end{proof}
%%%%%%%%%%%%%%%%%%%%%%%%%%%%%%%%%%%%%% 

Next, following Kato \cite{Ka66}, one introduces
\begin{align}
\begin{split}
& R(z)=R_0(z)-\ol{R_0(z)B^*}[I_{\cK}-K(z)]^{-1}AR_0(z), \lb{2.13} \\
& \hspace*{2.65cm} z\in\{\zeta\in\rho(H_0)\,|\, 1\in\rho(K(\zeta))\}.
\end{split}
\end{align}

%%%%%%%%%%%%%%%%%%%%%%%%%%%%%%%%%%%%%%
\begin{theorem} \lb{t2.3}
Assume Hypothesis \ref{h2.1} and suppose
$z\in\{\zeta\in\rho(H_0)\,|\, 1\in\rho(K(\zeta))\}$. Then, $R(z)$
defined in \eqref{2.13} defines a densely defined, closed, linear
operator $H$ in $\cH$ by
\begin{equation}
R(z)=(H-zI_{\cH})^{-1}. \lb{2.14}
\end{equation}
Moreover,
\begin{equation}
AR(z), BR(z)^* \in \cB(\cH,\cK) \lb{2.15}
\end{equation}
and
\begin{align}
R(z)&=R_0(z)-\ol{R(z)B^*}AR_0(z) \lb{2.16} \\
&=R_0(z)-\ol{R_0(z)B^*}AR(z).   \lb{2.17}
\end{align}
Finally, $H$ is an extension of
$(H_0+B^*A)|_{\dom(H_0)\cap\dom(B^*A)}$
$($the latter intersection domain may consist of $\{0\}$ only$)$,
\begin{equation}
H\supseteq (H_0+B^*A)|_{\dom(H_0)\cap\dom(B^*A)}. \lb{2.18}
\end{equation}
\end{theorem}
%%%%%%%%%%%%%%%%%%%%%%%%%%%%%%%%%%%%%
\begin{proof}
Suppose $z\in\{\zeta\in\rho(H_0)\,|\, 1\in\rho(K(\zeta))\}$. Since
by \eqref{2.13}
\begin{align}
AR(z)&=[I_{\cK}-K(z)]^{-1}AR_0(z), \lb{2.19} \\
BR(z)^*&=[I_{\cK}-K(z)^*]^{-1}BR_0(z)^*, \lb{2.20}
\end{align}
$R(z)f=0$ implies $AR(z)f=0$ and hence by \eqref{2.19} $AR_0(z)f=0$.
The latter implies $R_0(z)f=0$ by \eqref{2.13} and thus $f=0$.
Consequently,
\begin{equation}
\ker(R(z))=\{0\}. \lb{2.21}
\end{equation}
Similarly, \eqref{2.20} implies
\begin{equation}
\ker(R(z)^*)=\{0\} \, \text{ and hence } \, \ol{\ran(R(z))}=\cH.
\lb{2.22}
\end{equation}
Next, combining \eqref{2.13}, the resolvent equation for $H_0$,
\eqref{2.6}, \eqref{2.7}, \eqref{2.10}, and \eqref{2.11} proves the
resolvent equation
\begin{align}
\begin{split}
& R(z_1)-R(z_2)=(z_1-z_2)R(z_1)R(z_2), \lb{2.23} \\
& \hspace*{.6cm}
z_1, z_2 \in\{\zeta\in\rho(H_0)\,|\, 1\in\rho(K(\zeta))\}.
\end{split}
\end{align}
Thus, $R(z)$ is indeed the resolvent of a densely defined, closed,
linear operator $H$ in $\cH$ as claimed in connection with \eqref{2.14}.

By \eqref{2.19} and \eqref{2.20}, $AR(z)\in\cB(\cH,\cK)$ and
$[BR(\ol z)^*]^*=\ol{R(z)B^*}\in\cB(\cK,\cH)$, proving \eqref{2.15}. A
combination of \eqref{2.13}, \eqref{2.19}, and \eqref{2.20} then
proves \eqref{2.16} and \eqref{2.17}.

Finally, let $f\in\dom(H_0)\cap\dom(B^*A)$ and set $g=(H_0-z I_{\cH})f$.
Then $R_0(z)g=f$ and by \eqref{2.16}, $R(z)g-f=-R(z)B^*Af$.
Thus, $f\in\dom(H)$ and $(H-zI_{\cH})f=g+B^*Af
=(H_0+B^*A-zI_{\cH})f$, proving \eqref{2.18}.
\end{proof}
%%%%%%%%%%%%%%%%%%%%%%%%%%%%%%%%%%%%

%%%%%%%%%%%%%%%%%%%%%%%%%%%%%%%%%%%%
\begin{remark} \lb{r2.4}
$(i)$ Assume that $H_0$ is self-adjoint in $\cH$. Then $H$ is also
self-adjoint if
\begin{equation}
(Af,Bg)_{\cK}=(Bf,Ag)_{\cK} \, \text{ for all } \,
f,g\in\dom(A)\cap\dom(B). \lb{2.24}
\end{equation}
$(ii)$ The formalism is symmetric with respect to $H_0$ and $H$ in the
following sense: The densely defined operator $-AR(z)B^*$ has a
bounded extension to all of
$\cK$ for all $z\in\{\zeta\in\rho(H_0)\,|\, 1\in\rho(K(\zeta))\}$, in
particular,
\begin{equation}
I_{\cK}-\ol{AR(z)B^*}=[I_{\cK}-K(z)]^{-1}, \quad
z\in\{\zeta\in\rho(H_0)\,|\, 1\in\rho(K(\zeta))\}. \lb{2.25}
\end{equation}
Moreover,
\begin{align}
\begin{split}
R_0(z)&=R(z)+\ol{R(z)B^*}[I_{\cK}-\ol{AR(z)B^*}]^{-1}AR(z), \\
& \hspace*{2.25cm}
z\in\{\zeta\in\rho(H_0)\,|\, 1\in\rho(K(\zeta))\} \lb{2.26}
\end{split}
\end{align}
and
\begin{equation}
H_0\supseteq (H-B^*A)|_{\dom(H)\cap\dom(B^*A)}.  \lb{2.27}
\end{equation}
$(iii)$ The basic hypotheses \eqref{2.2} which amount to
\begin{equation}
AR_0(z)\in\cB(\cH,\cK), \quad \ol{R_0(z)B^*}=[B(H_0^*-\ol z)^{-1}]^* \in
\cB(\cK,\cH), \quad z\in\rho(H_0)  \lb{2.27a}
\end{equation}
(cf.\ \eqref{2.5}) are more general than a quadratic form perturbation
approach which would result in conditions of the form
\begin{equation}
AR_0(z)^{1/2}\in\cB(\cH,\cK), \quad \ol{R_0(z)^{1/2}B^*}=[B(H_0^*-\ol
z)^{-1/2}]^* \in \cB(\cK,\cH), \quad z\in\rho(H_0),  \lb{2.27b}
\end{equation}
or even an operator perturbation approach which would involve conditions
of the form
\begin{equation}
[B^* A] R_0(z)\in\cB(\cH), \quad z\in\rho(H_0).  \lb{2.27c}
\end{equation}
${}$ \hfill $\diamond$
\end{remark}
%%%%%%%%%%%%%%%%%%%%%%%%%%%%%%%%%%%%%

%%%%%%%%%%%%%%%%%%%%%%%%%%%%%%%%%%%%%
%%%%%%%%%%%%%%%%%%%%%%%%%%%%%%%%%%%%% 
\section{A General Birman--Schwinger Principle} \lb{s3}
%%%%%%%%%%%%%%%%%%%%%%%%%%%%%%%%%%%%%
%%%%%%%%%%%%%%%%%%%%%%%%%%%%%%%%%%%%%

The principal result in this section represents an abstract version
of (a variant of) the Birman--Schwinger principle due to Birman \cite{Bi61}
and Schwinger \cite{Sc61} (cf.\ also \cite{BS91}, \cite{GH87},
\cite{Kl82}, \cite{KS80}, \cite{Ne83}, \cite{Ra80}, \cite{Se74}, and
\cite{Si77a}).

We need to strengthen our hypotheses a bit and hence introduce the following
assumption:

%%%%%%%%%%%%%%%%%%%%%%%%%%%%%%%%%%%%%
\begin{hypothesis} \lb{h2.5}
In addition to Hypothesis \ref{h2.1} we suppose the condition: \\
$(iv)$ $K(z)\in\cB_\infty(\cK)$ for all $z\in\rho(H_0)$.
\end{hypothesis}
%%%%%%%%%%%%%%%%%%%%%%%%%%%%%%%%%%%%%

Since by \eqref{2.25}
\begin{align}
-\ol{AR(z)B^*}&=[I_{\cK}-K(z)]^{-1}K(z) \lb{2.28} \\
&=-I_{\cK}+[I_{\cK}-K(z)]^{-1}, \lb{2.28a}
\end{align}
Hypothesis \ref{h2.5} implies that $-\ol{AR(z)B^*}$ extends to a
compact operator in $\cK$ as long as the right-hand side of
\eqref{2.28a} exists.

The following general result is due to Konno and Kuroda \cite{KK66}
in the case where $H_0$ is self-adjoint. (The more general case
presented here requires no modifications but we present a proof for
completeness.)
%%%%%%%%%%%%%%%%%%%%%%%%%%%%%%%%%%%% 
\begin{theorem}[\cite{KK66}] \lb{t2.6}
Assume Hypothesis \ref{h2.5} and let $\lambda_0\in\rho(H_0)$. Then,
\begin{equation}
Hf=\lambda_0 f, \quad 0\neq f\in\dom(H) \, \text{ implies } \,
K(\lambda_0)g=g \lb{2.29}
\end{equation}
where, for fixed $z_0\in\{\zeta\in\rho(H_0)\,|\,
1\in\rho(K(\zeta))\}$, $z_0\neq \lambda_0$,
\begin{align}
0 \neq g &= [I_{\cK}-K(z_0)]^{-1}AR_0(z_0)f \lb{2.30} \\
&= (\lambda_0-z_0)^{-1}Af. \lb{2.30a}
\end{align}
Conversely,
\begin{equation}
K(\lambda_0)g=g, \quad 0\neq g\in\cK \, \text{ implies } \,
Hf=\lambda_0 f, \lb{2.31}
\end{equation}
where
\begin{equation}
0\neq f=-\ol{R_0(\lambda_0)B^*}g\in\dom(H). \lb{2.32}
\end{equation}
Moreover,
\begin{equation}
\dim(\ker(H-\lambda_0I_{\cH}))
=\dim(\ker(I_{\cK}-K(\lambda_0)))<\infty.
\lb{2.33}
\end{equation}
In particular, let $z\in\rho(H_0)$, then
\begin{equation}
\text{$z\in\rho(H)$ if and only if\;
$1\in\rho(K(z))$.} \lb{2.33a}
\end{equation}
\end{theorem}
%%%%%%%%%%%%%%%%%%%%%%%%%%%%%%%%%%%
\begin{proof}
$Hf=\lambda_0f$, $0\neq f\in\dom(H)$, is equivalent to
$f=(\lambda_0-z_0)R(z_0)f$ and applying \eqref{2.13} one obtains
after a simple rearrangement that
\begin{equation}
(H_0-\lambda_0I_{\cH})R_0(z_0)f=-(\lambda_0-z_0)\ol{R_0(z_0)B^*}
[I_{\cK}-K(z_0)]^{-1}AR_0(z_0)f. \lb{2.34}
\end{equation}
Next, define $g=[I_{\cK}-K(z_0)]^{-1}AR_0(z_0)f$. Then $g\neq 0$
since otherwise
\begin{equation}
(H_0-\lambda_0 I_{\cH})R_0(z_0)f=0, \quad 0 \neq R_0(z_0)f\in \dom(H_0),
\, \text{ and hence } \, \lambda_0\in\sigma(H_0), \lb{2.35}
\end{equation}
would contradict our hypothesis $\lambda_0\in\rho(H_0)$. Applying
$[I_{\cK}-K(z_0)]^{-1}AR_0(\lambda_0)$ to \eqref{2.34} then yields
\begin{align}
&[I_{\cK}-K(z_0)]^{-1}AR_0(\lambda_0)(H_0-\lambda_0 I_{\cH})R_0(z_0)f=
[I_{\cK}-K(z_0)]^{-1}AR_0(z_0)f=g \no \\
&\quad =-(\lambda_0-z_0)[I_{\cK}-K(z_0)]^{-1}AR_0(\lambda_0)
\ol{R_0(z_0)B^*}[I_{\cK}-K(z_0)]^{-1}AR_0(z_0)f \no \\
&\quad =-(\lambda_0-z_0)[I_{\cK}-K(z_0)]^{-1}AR_0(\lambda_0)
\ol{R_0(z_0)B^*}g.  \lb{2.36}
\end{align}
Thus, based on \eqref{2.10}, one infers
\begin{align}
g&=-(\lambda_0-z_0)[I_{\cK}-K(z_0)]^{-1}AR_0(\lambda_0)
\ol{R_0(z_0)B^*}g \no \\
&=[I_{\cK}-K(z_0)]^{-1}[K(\lambda_0)-K(z_0)]g \no \\
&=g-[I_{\cK}-K(z_0)]^{-1}[I_{\cK}-K(\lambda_0)]g \lb{2.37}
\end{align}
and hence $K(\lambda_0)g=g$, proving \eqref{2.29}. Since
$f=(\lambda_0-z_0)R(z_0)f$, using \eqref{2.19} one computes
\begin{align}
Af&=(\lambda_0-z_0)AR(z_0)f \no \\
&=(\lambda_0-z_0)[I_{\cK}-K(z_0)]^{-1}AR_0(z_0)f \no \\
&=(\lambda_0-z_0)g,  \lb{2.37a}
\end{align}
proving \eqref{2.30a}.

Conversely, suppose $K(\lambda_0)g=g$, $0\neq g\in\cK$ and define
$f=-\ol{R_0(\lambda_0)B^*}g$. Then a simple computation using
\eqref{2.10} shows
\begin{align}
g&=g-[I_{\cK}-K(z_0)]^{-1}[I_{\cK}-K(\lambda_0)]g \no \\
&=[I_{\cK}-K(z_0)]^{-1}[K(\lambda_0)-K(z_0)]g \no \\
&=(\lambda_0-z_0)[I_{\cK}-K(z_0)]^{-1}AR_0(z_0)f. \lb{2.38}
\end{align}
Thus, $f\neq0$ since $f=0$ would imply the contradiction $g=0$. Next,
inserting the definition of $f$ into \eqref{2.38} yields
\begin{align}
g&=(\lambda_0-z_0)[I_{\cK}-K(z_0)]^{-1}AR_0(z_0)f \no \\
&=-(\lambda_0-z_0)[I_{\cK}-K(z_0)]^{-1}AR_0(z_0)
\ol{R_0(\lambda_0)B^*}g. \lb{2.39}
\end{align}
Applying $\ol{R_0(z_0)B^*}$ to \eqref{2.39} and taking into account
\begin{align}
\ol{R_0(z_0)B^*}g&=\ol{[R_0(\lambda_0)
-(\lambda_0-z_0)R_0(z_0)R_0(\lambda_0)]B^*}g \no \\
&=-f+(\lambda_0-z_0)R_0(z_0)f,  \lb{2.40}
\end{align}
a combination of \eqref{2.40} and \eqref{2.13} yields that
\begin{align}
& -f-(z_0-\lambda_0)R_0(z_0)f=(\lambda_0-z_0)\ol{R_0(z_0)B^*}
[I_{\cK}-K(z_0)]^{-1}AR_0(z_0)f \no \\
& \quad = (\lambda_0-z_0)[R_0(z_0)-R(z_0)]f.   \lb{2.41}
\end{align}
The latter is equivalent to $(\lambda_0-z_0)(H-z_0I_{\cH})^{-1}f=f$.
Thus, $f\in\dom(H)$ and $Hf=\lambda_0f$, proving \eqref{2.31}.

Since $K(\lambda_0)\in\cB_\infty(\cK)$, the eigenspace of
$K(\lambda_0)$ corresponding to the eigenvalue $1$ is
finite-dimensional. The previous considerations established a
one-to-one correspondence between the geometric eigenspace of
$K(\lambda_0)$ corresponding to the eigenvalue $1$ and the geometric
eigenspace of $H$ corresponding to the eigenvalue $\lambda_0$. This
proves \eqref{2.33}.

Finally, \eqref{2.33}, \eqref{2.13}, and \eqref{2.25} prove
\eqref{2.33a}.
\end{proof}
%%%%%%%%%%%%%%%%%%%%%%%%%%%%%%%%%%%%%

%%%%%%%%%%%%%%%%%%%%%%%%%%%%%%%%%%%%%
\begin{remark}
It is possible to avoid the compactness assumption in Hypothesis \ref{h2.5}
in Theorem \ref{t2.6} provided that \eqref{2.33} is replaced by the
statement
\begin{equation}\label{2.33NEW}
\text{the subspaces $\ker(H-\lambda_0I_{\cH})$ and
$\ker(I_{\cK}-K(\lambda_0))$ are isomorphic.}
\end{equation}
(Of course, \eqref{2.33} follows from \eqref{2.33NEW} provided
$\ker(I_{\cK}-K(\lambda_0))$ is finite-dimensional, which in turn follows
from Hypothesis \ref{h2.5}). Indeed, by formula \eqref{2.19}, we have
$AR(z_0)=[I_{\cK}-K(z_0)]^{-1}AR_0(z_0)$. By formula \eqref{2.30},
if $f\neq 0$, then $g=AR(z_0)f\neq 0$, and thus the operator
\begin{equation}
AR(z_0)=[I_{\cK}-K(z_0)]^{-1}AR_0(z_0) \colon
\ker(H-\lambda_0I)\to\ker(K(\lambda_0)-I)
\end{equation}
is injective. By formula \eqref{2.39} this operator is also surjective,
since each $g\in\ker(K(\lambda_0)-I)$ belongs to its range,
\begin{equation}
g=(\lambda_0-z_0)[I_{\cK}-K(z_0)]^{-1}AR_0(z_0)f=AR(z_0)f,
\end{equation}
where $f\in\ker(H-\lambda_0I)$. \hfill $\diamond$
\end{remark}
%%%%%%%%%%%%%%%%%%%%%%%%%%%%%%%%%%%%%%

%%%%%%%%%%%%%%%%%%%%%%%%%%%%%%%%%%%%%%
%%%%%%%%%%%%%%%%%%%%%%%%%%%%%%%%%%%%%% 
\section{Essential Spectra and a Local Weinstein--Aronszajn Formula}
\lb{s4}
%%%%%%%%%%%%%%%%%%%%%%%%%%%%%%%%%%%%%%
%%%%%%%%%%%%%%%%%%%%%%%%%%%%%%%%%%%%%%

In this section, we closely follow Howland \cite{Ho70} and prove
a result which demonstrates the invariance of the essential spectrum.
However, since we will extend Howland's result to the
non-self-adjoint case, this requires further explanation. Moreover, we
will also re-derive Howland's local Weinstein--Aronszajn formula.

%%%%%%%%%%%%%%%%%%%%%%%%%%%%%%%%%%%%%%
\begin{definition} \lb{d2.7}
Let $\Omega\subseteq\bbC$ be open and connected. Suppose
$\{L(z)\}_{z\in\Omega}$ is a family of compact operators in $\cK$, which is
analytic on $\Omega$ except for isolated singularities. Following Howland
we call $\{L(z)\}_{z\in\Omega}$ {\it completely meromorphic} on $\Omega$ if
$L$ is meromorphic on $\Omega$ and the principal part of
$L$ at each of its poles is of finite rank.
\end{definition}
%%%%%%%%%%%%%%%%%%%%%%%%%%%%%%%%%%%%%%

We start with an auxiliary result due to Steinberg
\cite{St68} with a modification by Howland \cite{Ho70}.

%%%%%%%%%%%%%%%%%%%%%%%%%%%%%%%%%%%%%%
\begin{lemma}[\cite{Ho70}, \cite{St68}] \lb{l2.8}
Let $\{L(z)\}_{z\in\Omega}$ be an analytic $($resp., completely
meromor-\ phic$)$ family of compact operators in $\cK$ on an open connected set
$\Omega\subseteq\bbC$. Then for each $z_0\in\Omega$ there is a
neighborhood $U(z_0)$ of $z_0$, and an analytic $\cB(\cK)$-valued function
$M$ on $U(z_0)$, such that $M(z)^{-1}\in \cB(\cK)$ for all $z\in U(z_0)$
and
\begin{equation}
M(z)[I_{\cK}-L(z)]=I_{\cK}-F(z), \quad z\in U(z_0), \lb{2.44}
\end{equation}
where $F$ is analytic $($resp., meromorphic$)$ on $U(z_0)$
with $F(z)$ of finite rank $($except at poles$)$ for all $z\in U(z_0)$.
\end{lemma}
%%%%%%%%%%%%%%%%%%%%%%%%%%%%%%%%%%%%%%

The next auxiliary result is due to Ribaric and Vidav \cite{RV69}.

%%%%%%%%%%%%%%%%%%%%%%%%%%%%%%%%%%%%%%
\begin{lemma}[\cite{RV69}] \lb{l2.9}
Let $\{L(z)\}_{z\in\Omega}$ be a completely meromorphic family in $\cK$ on
an open connected set $\Omega\subseteq\bbC$. Then either \\
$(i)$ $I_{\cK}-L(z)$ is not boundedly invertible for all $z\in\Omega$,
\\
or \\
$(ii)$ $\{[I_{\cK}-L(z)]^{-1}-I_{\cK}\}_{z\in\Omega}$ is completely
meromorphic on $\Omega$.
\end{lemma}
%%%%%%%%%%%%%%%%%%%%%%%%%%%%%%%%%%%%%

Moreover, we state the following result due to Howland \cite{Ho71}.

%%%%%%%%%%%%%%%%%%%%%%%%%%%%%%%%%%%%% 
\begin{lemma}[\cite{Ho71}] \lb{l2.10}
Let $\{L(z)\}_{z\in\Omega}$ be an analytic $($resp., meromorphic$)$ family
in $\cK$ on an open connected set $\Omega\subseteq\bbC$ and suppose
that $L(z)$ has finite rank for each $z\in\Omega$ $($except at
poles$)$. Then the following assertions hold: \\
$(i)$ \, The rank of $L(z)$ is constant for all $z\in\Omega$, except for
isolated points where it decreases. \\
$(ii)$ \, $\Delta(z)=\det(I_{\cK}-L(z))$ and \,$\tr(L(z))$
are analytic $($resp., meromorphic$)$ for all $z\in\Omega$. \\
$(iii)$ Whenever $\Delta(z)\neq 0$,
\begin{equation}
\Delta'(z)/\Delta(z)=-\tr([I_{\cK}-L(z)]^{-1}L'(z)),
\quad z\in\Omega. \lb{2.44a}
\end{equation}
\end{lemma}
%%%%%%%%%%%%%%%%%%%%%%%%%%%%%%%%%%%%%

We note that it can of course happen that $\Delta$ vanishes
identically on $\Omega$.

Next, we introduce the multiplicity function $m(\cdot,T)$ on
$\bbC$ associated with a closed, densely defined, linear operator
$T$ in $\cH$ as follows. Suppose $\lambda_0\in\bbC$ is an isolated
point in $\sigma(T)$ and introduce the Riesz projection
$P(\lambda_0,T)$ of $T$ corresponding to $\lambda_0$ by
\begin{equation}
P(\lambda_0,T)=-\f{1}{2\pi i}\oint_{C(\lambda_0; \varepsilon) }
d\zeta \, (T-\zeta I_{\cH})^{-1}, \lb{2.45}
\end{equation}
where $C(\lambda_0; \varepsilon) $ is a counterclockwise oriented
circle centered at $\lambda_0$ with sufficiently small radius
$\varepsilon>0$ (excluding the rest of $\sigma(T)$). Then
$m(z,T)$, $z\in\bbC$, is defined by
\begin{equation}
m(z,T)=\begin{cases} 0, & \text{if $z\in\rho(T)$,} \\
\dim(\ran(P(z,T))), & \text{if $z$ is an isolated eigenvalue of
$T$} \\
& \text{of finite algebraic multiplicity,} \\
+\infty, & \text{otherwise.} \end{cases} \lb{2.46}
\end{equation}
(Here ``isolated eigenvalue" of $T$ means an isolated point in the spectrum 
of $T$ which is an eigenvalue of $T$.)  
We note that the dimension of the Riesz projection in \eqref{2.45} is
finite if and only if $\lambda_0$ is an isolated eigenvalue of $T$ of
finite algebraic multiplicity (cf.\ \cite[p.\ 181]{Ka80}). In analogy
to the self-adjoint case (but deviating from most definitions in the
non-self-adjoint case, see \cite[Sect.\ I.4, Ch.\ IX]{EE89}) we now
introduce the set 
\begin{align}
\begin{split}
\wti \sigma_{\rm e}(T)&=\{\lambda\in\bbC\, | \, \lambda\in\sigma(T),
\, \text{$\lambda$ is not an isolated eigenvalue of $T$}   \\
& \hspace*{4.2cm} \text{of finite algebraic multiplicity}\}. \lb{2.47}
\end{split}
\end{align}
Of course, $\wti \sigma_{\rm e}(T)$ coincides with the essential
spectrum of $T$ if $T$ is self-adjoint in $\cH$. In the
non-self-adjoint case at hand, the set $\wti \sigma_{\rm e}(T)$ is
most natural in our study of $H_0$ and $H$ as will
subsequently be shown. It will also be convenient to introduce the
complement of $\wti \sigma_{\rm e}(T)$ in $\bbC$,
\begin{align}
\wti \rho (T)&=\bbC\backslash\wti\sigma_{\rm e}(T)    \no \\
&=\rho(T)\cup\{\lambda\in\bbC\,|\, \text{$\lambda$ is an isolated 
eigenvalue of $T$ of}     \lb{2.48A}  \\
&\hspace*{3.9cm} \text{finite algebraic multiplicity}\}.  \no
\end{align}

If $\lambda_0\in\bbC$ is an isolated eigenvalue of $T$
of finite algebraic multiplicity, then the singularity structure of
the resolvent of $T$ near $\lambda_0$ is of the type
\begin{align}
(T-zI_{\cH})^{-1}&=(\lambda_0-z)^{-1}P(\lambda_0,T)
+\sum_{k=1}^{\mu(\lambda_0,T)} (\lambda_0-z)^{-k-1}
(-1)^kD(\lambda_0,T)^k \no \\
& \quad +\sum_{k=0}^\infty (\lambda_0-z)^k (-1)^k
S(\lambda_0,T)^{k+1} \lb{2.48a}
\end{align}
for $z$ in a sufficiently small neighborhood of $\lambda_0$. Here
\begin{align}
D(\lambda_0,T)&=(T-\lambda_0 I_{\cH})P(\lambda_0,T)=\f{1}{2\pi i}
\oint_{C(\lambda_0; \varepsilon) } d\zeta \,
(\lambda_0-\zeta)(T-\zeta I_{\cH})^{-1} \in \cB(\cH),
\lb{2.48b} \\
S(\lambda_0,T)&=-\f{1}{2\pi i} \oint_{C(\lambda_0; \varepsilon) }
d\zeta \, (\lambda_0-\zeta)^{-1}(T-\zeta I_{\cH})^{-1} \in
\cB(\cH), \lb{2.48c}
\end{align}
and $D(\lambda_0,T)$ is nilpotent with its range contained in that of
$P(\lambda_0,T)$,
\begin{equation}
D(\lambda_0,T)=P(\lambda_0,T)D(\lambda_0,T)
=D(\lambda_0,T)P(\lambda_0,T). \lb{2.48d}
\end{equation}
Moreover,
\begin{align}
\begin{split}
&S(\lambda_0,T)T\subseteq TS(\lambda_0,T), \quad
(T-\lambda_0I_{\cH})S(\lambda_0,T)
= I_{\cH}-P(\lambda_0,T), \\
&S(\lambda_0,T)P(\lambda_0,T)=
P(\lambda_0,T)S(\lambda_0,T)=0. \lb{2.48e}
\end{split}
\end{align}
Finally,
\begin{align}
&\mu(\lambda_0,T) + 1 \leq m(\lambda_0,T)=\dim(\ran(P(\lambda_0,T))),
\lb{2.48f} \\
&\tr(P(\lambda_0,T))=m(\lambda_0,T), \quad D(\lambda_0,T)^{m(\lambda_0,T)} = 0. \lb{2.48g}
\end{align}

Next, we need one more notation: Let
$\Omega\subseteq\bbC$ be open and connected, and let
$f\colon\Omega\to\bbC\cup\{\infty\}$ be meromorphic and not
identically vanishing on $\Omega$. The multiplicity function
$m(z;f)$, $z\in\Omega$, is then defined by
\begin{align}
m(z;f)&=\begin{cases} k, & \text{if $z$ is a zero of $f$ of order
$k$,} \\
-k, & \text{if $z$ is a pole of order $k$,} \\
0, & \text{otherwise.} \end{cases} \lb{2.49} \\
&= \f{1}{2\pi i}\oint_{C(z; \varepsilon) } d\zeta \,
\f{f'(\zeta)}{f(\zeta)}, \quad z\in\Omega \lb{2.49a}
\intertext{for $\varepsilon>0$
sufficiently small. If $f$ vanishes identically on $\Omega$, one
defines} m(z;f)&=+\infty,  \quad z\in\Omega. \lb{2.50}
\end{align}
Here the circle $C(z; \varepsilon) $ is chosen sufficiently small such
that $C(z; \varepsilon) $ contains no other singularities or zeros of
$f$ except, possibly, $z$.

The following result is due to Howland in the case where $H_0$ and
$H$ are self-adjoint. We will closely follow his strategy of proof
and present detailed arguments in the more general situation considered
here.

%%%%%%%%%%%%%%%%%%%%%%%%%%%%%%%%%%%%%%
\begin{theorem} \lb{t2.11}
In addition to Hypothesis \ref{h2.5} assume that each connected component of $\widetilde \rho(H_0)$ 
contains a point of $\rho(H)$, and each connected component of $\widetilde \rho(H)$ 
contains a point of $\rho(H_0)$. Then\footnote{The additional hypotheses on the connected components of 
$\widetilde \rho(H_0)$ and $\widetilde \rho(H)$ were originally missed in our 2005 paper published in Russ. J. Math. Phys. {\bf 12}, 443--471 (2005). The necessity of these additional hypothesis for \eqref{2.51} to hold is discussed in detail in \cite[Sect.~6]{BEG20}. \label{F}},
\begin{equation}
\wti\sigma_{\rm e}(H)=\wti\sigma_{\rm e}(H_0). \lb{2.51}
\end{equation}
In addition, let $\lambda_0\in
\bbC\backslash\wti\sigma_{\rm e}(H_0)$. Then
there exists a neighborhood $U(\lambda_0)$ of $\lambda_0$ and a
function $\Delta(\cdot)$ meromorphic on $U(\lambda_0)$, which does
not vanish identically, such that the local Weinstein--Aronszajn
formula
\begin{equation}
m(z,H)=m(z,H_0)+m(z;\Delta), \quad z\in U(\lambda_0) \lb{2.52}
\end{equation}
holds.
\end{theorem}
%%%%%%%%%%%%%%%%%%%%%%%%%%%%%%%%%%%%
\begin{proof}
By \eqref{2.10}, $K(\cdot)$ is analytic on $\rho(H_0)$ and
\begin{equation}
K'(z)=-AR_0(z)[BR_0(z)^*]^*, \quad z\in\rho(H_0). \lb{2.63}
\end{equation}
Let $z_0\in\wti \rho (H_0)$, then by \eqref{2.48a},
\begin{align}
\begin{split} 
& R_0(z)=(z_0-z)^{-1}P_0 +
\sum_{k=1}^{\mu_0} (z_0-z)^{-k-1}(-1)^k D_0^k + G_0(z) \lb{2.53} \\
& \hspace*{3.85cm} \text{for $0 < |z - z_0|$ sufficiently small,} 
\end{split}
\end{align}
where $G_0(\cdot)$ is analytic in a neighborhood of $z_0$. Since
\begin{equation}
\ran(D_0)\subseteq\ran(P_0)\subset\dom(H_0)\subseteq \dom(A),
\lb{2.54}
\end{equation}
$AP_0B^*$, $AD_0B^*$, and $AG_0(z)B^*$ have compact extensions from
$\dom(B^*)$ to $\cK$, and the extensions of $AP_0B^*$ and $AD_0B^*$
are given by the finite-rank operators $AP_0[BP_0^*]^*$ and
$\ol{AP_0D_0P_0B^*}$, respectively. Moreover, 
the compact extension of $AG_0(z)B^*$ is analytic near $z_0$. Consequently,
$K(\cdot)$ is completely meromorphic on the connected component of $\wti \rho (H_0)$.

Similarly, by \eqref{2.28a}, the fact that every connected component of $\widetilde \rho(H_0)$ 
contains a point of $\rho(H)$, and Lemma \ref{l2.9}, $-\ol{AR(z)B^*}$
is completely meromorphic on $\wti \rho (H_0)$. Moreover, by
\eqref{2.28a}, any singularity $z_0$ of $-\ol{AR(z)B^*}$ is an
isolated point of $\sigma(H)$. Since $R_0(z)$, $AR_0(z)$, and
$BR_0(z)$ all have finite-rank principal parts at their poles,
\eqref{2.13} and \eqref{2.28a} show that $R(z)$ also has a
finite-rank principal part at $z_0$. The latter implies that $z_0$ is an
eigenvalue of $H$ of finite algebraic multiplicity. Thus,
$\wti \rho (H_0)\subseteq\wti \rho (H)$. Since by Remark \ref{r2.4}\,$(ii)$ and the assumptions 
in Theorem \ref{t2.11} this formalism is symmetric with respect to $H_0$ and
$H$, one also obtains $\wti \rho (H_0)\supseteq\wti \rho (H)$, and
hence \eqref{2.51}. 

Next, by Lemma \ref{l2.8}, let $U_0$ be a neighborhood of $\lambda_0$
such that
\begin{equation}
M(z)[I_{\cK}-K(z)]=I_{\cK}-F(z), \quad \lb{2.67}
\end{equation}
with $M$ analytic and boundedly invertible on $U_0$ and some
$F$ meromorphic and of finite rank on $U_0$. One defines
\begin{equation}
\Delta(z)=\det(I_{\cK}-F(z)), \quad z\in U_0. \lb{2.68}
\end{equation}
Since by Lemma \ref{l2.9}, $[I_{\cK}-K(z)]^{-1}$ is meromorphic
and $M(z)$ is boundedly invertible for all $z\in U_0$,
$[I_{\cK}-F(z)]^{-1}$ is also meromorphic on $U_0$, and hence,
$\Delta(\cdot)$ is not identically zero on $U_0$. By Lemma
\ref{l2.10}\,(iii) and cyclicity of the trace (i.e., $\tr(ST)=\tr(TS)$ for
$S$ and $T$ bounded operators such that $ST$ and $TS$ lie in the trace
class, cf.\ \cite[Corollary 3.8]{Si79}),
\begin{align}
\Delta'(z)/\Delta(z)&=-\tr([I_{\cK}-F(z)]^{-1}F'(z)) \no \\
&=\tr([I_{\cK}-K(z)]^{-1}M(z)^{-1}M'(z)[I_{\cK}-K(z)]-
[I_{\cK}-K(z)]^{-1}K'(z)) \no \\
&=\tr(M(z)^{-1}M'(z)-K'(z)[I_{\cK}-K(z)]^{-1}). \lb{2.69}
\end{align}
Let $z_0\in U_0$ and $C(z_0; \varepsilon)$ be a clockwise oriented
circle centered at $z_0$ with sufficiently small radius
$\varepsilon$ (excluding all singularities of $[I_{\cK}-F(z)]^{-1}$,
except, possibly, $z_0$) contained in $U_0$. Then,
\begin{align}
m(z_0;\Delta)&=\f{1}{2\pi i} \oint_{C(z_0; \varepsilon) } d\zeta \,
\f{\Delta'(\zeta)}{\Delta(\zeta)} \no \\
&= \f{1}{2\pi i} \oint_{C(z_0; \varepsilon) } d\zeta \,
\tr(M(\zeta)^{-1}M'(\zeta)-K'(\zeta)[I_{\cK}-K(\zeta)]^{-1}).
\lb{2.70}
\end{align}
Since $M$ is analytic and boundedly invertible on $U_0$, an interchange of
the trace and the integral, using
\begin{equation}
\oint_{C(z_0; \varepsilon) } d\zeta \,
M(\zeta)^{-1}M'(\zeta) =0
\end{equation}
and \eqref{2.63}, then yields
\begin{align}
m(z_0;\Delta)&= \f{1}{2\pi i} \tr\bigg(\oint_{C(z_0; \varepsilon) }
d\zeta \, AR_0(\zeta)[BR_0(\zeta)^*]^*[I_{\cK}-K(\zeta)]^{-1}
\bigg)
\no \\
&=\f{1}{2\pi i} \tr\bigg(\oint_{C(z_0; \varepsilon) } d\zeta \,
AR_0(\zeta)[BR(\zeta)^*]^*\bigg). \lb{2.71}
\end{align}
Next, for $\varepsilon>0$ sufficiently small, one infers from
\cite[p.\ 178]{Ka80} (cf.\ \eqref{2.48g}) that
\begin{align}
m(z_0,H)-m(z_0,H_0)&=-\f{1}{2\pi i} \tr\bigg(
\oint_{C(z_0; \varepsilon) } d\zeta \, [R(\zeta)-R_0(\zeta)]\bigg)
\no
\\
&= \f{1}{2\pi i} \tr\bigg( \oint_{C(z_0; \varepsilon) } d\zeta \,
[BR_0(\zeta)^*]^*
[I_{\cK}-K(\zeta)]^{-1} AR_0(\zeta) \bigg) \no \\
&= \f{1}{2\pi i} \tr\bigg( \oint_{C(z_0; \varepsilon) } d\zeta \,
[BR(\zeta)^*]^* AR_0(\zeta) \bigg). \lb{2.72}
\end{align}
At this point we cannot simply change back the order of the trace
and the integral and use the cyclicity of the trace to prove equality of
\eqref{2.71} and
\eqref{2.72} since now the integrand is not necessarily trace class.
But one can prove the equality of \eqref{2.71} and \eqref{2.72} directly
as follows. Writing (cf.\ \eqref{2.48a}),
\begin{align}
AR_0(z) &= (z_0-z)^{-1}\wti
P_0+\sum_{k=1}^{\mu_0}(z_0-z)^{-k-1}(-1)^k \wti D_0^k +
\sum_{k=0}^{\infty} (z_0-z)^k (-1)^k \wti S_0^{k+1}, \lb{2.73}
\\  [BR(z)^*]^* &= (z_0-z)^{-1}\wti Q_0
+\sum_{k=1}^{\nu_0}(z_0-z)^{-k-1}(-1)^k \wti E_0^k +
\sum_{k=0}^{\infty} (z_0-z)^k (-1)^k \wti T_0^{k+1}, \lb{2.74}
\end{align}
one obtains
\begin{align}
\text{\rm res}_{z=z_0}(AR_0(z)[BR(z)^*]^*)&=\wti P_0\wti T_0+\wti
S_0\wti Q_0+\sum_{k=1}^{\mu_0} \wti D_0^k\wti T_0^{k+1} +
\sum_{k=1}^{\nu_0} \wti S_0^{k+1}\wti E_0^{k}, \lb{2.75} \\
\text{\rm res}_{z=z_0}([BR(z)^*]^*AR_0(z))&=\wti T_0\wti P_0+\wti
Q_0\wti S_0+\sum_{k=1}^{\mu_0} \wti T_0^{k+1}\wti D_0^k +
\sum_{k=1}^{\nu_0} \wti E_0^{k}\wti S_0^{k+1}. \lb{2.76}
\end{align}
Using the cyclicity of the trace and Cauchy's theorem then
proves equality of \eqref{2.71} and \eqref{2.72} and hence
\eqref{2.52}.
\end{proof}
%%%%%%%%%%%%%%%%%%%%%%%%%%%%%%%%%%%%

%%%%%%%%%%%%%%%%%%%%%%%%%%%%%%%%%%%%
\begin{remark} \lb{r2.12}
Let $H_0$ be as in Hypothesis \ref{h2.1} and assume the hypotheses on 
$\widetilde \rho(H_0)$ and $\widetilde \rho(H)$ as in Theorem \ref{t2.11} 
(cf.\ footnote~\ref{F} on p.~11). \\
$(i)$ Let $V\in\cB_{\infty}(\cH)$ and define $H=H_0+V$,
$\dom(H)=\dom(H_0)$. Then \eqref{2.52} holds identifying $A=V$,
$B=I_{\cH}$, and $K(z)=VR_0(z)$ in connection with \eqref{2.13}. \\
$(ii)$ Let $V$ be of finite-rank and define $H=H_0+V$,
$\dom(H)=\dom(H_0)$. Then \eqref{2.52} holds on $\wti \rho (H_0)$
with $\Delta(z)=\det(I_{\cK}-K(z))$, $K(z)=VR_0(z)$, $z\in\rho(H_0)$, and
$U(\lambda_0)=\wti \rho (H_0)$. \hfill $\diamond$
\end{remark}
%%%%%%%%%%%%%%%%%%%%%%%%%%%%%%%%%%%%%

With the exception of the case discussed in Remark \ref{r2.12}\,(ii),
Theorem \ref{t2.11} has the drawback that it yields a
Weinstein--Aronszajn-type formula only locally on $U(\lambda_0)$.
However, by the same token, the great generality of this formalism,
basically assuming only compactness of $K(\cdot)$, must be
emphasized. In the following section we will present Howland's
global Aronszajn--Weinstein formula.

%%%%%%%%%%%%%%%%%%%%%%%%%%%%%%%%%%%%%%
%%%%%%%%%%%%%%%%%%%%%%%%%%%%%%%%%%%%%%
\section{A Global Weinstein--Aronszajn Formula} \lb{s5}
%%%%%%%%%%%%%%%%%%%%%%%%%%%%%%%%%%%%%%
%%%%%%%%%%%%%%%%%%%%%%%%%%%%%%%%%%%%%%

To this end we introduce a new hypothesis on $K$:

%%%%%%%%%%%%%%%%%%%%%%%%%%%%%%%%%%%%%%
\begin{hypothesis} \lb{h2.13}
In addition to the hypotheses in Theorem \ref{t2.11} (cf.\ footnote~\ref{F} on p.~11)  
we suppose the condition:  \\
$(v)$ For some $p\in\bbN$, $K(z)\in\cB_p(\cK)$ for all
$z\in\rho(H_0)$.
\end{hypothesis}
%%%%%%%%%%%%%%%%%%%%%%%%%%%%%%%%%%%%%%

We denote by $\|\cdot\|_p$ the norm in $\cB_p(\cK)$ and by
${\det}_p(\cdot)$ the regularized determinant of operators of the
type $I_{\cK}-L$, $L\in\cB_p(\cK)$ (cf.\ \cite{GGK96}, \cite{GGK97},
\cite[Chs.\ IX--XI]{GGK00}, \cite[Sect.\ 4.2]{GK69}, \cite{Si77},
\cite[Ch.\ 9]{Si79}).

We start by recalling the following result (cf.\ \cite[p.\
162--163]{GK69}, \cite[p.\ 107]{Si79}).

%%%%%%%%%%%%%%%%%%%%%%%%%%%%%%%%%%%%%%
\begin{lemma} \lb{l2.14}
Let $p\in\bbN$ and assume that $\{L(z)\}_{z\in\Omega}\in\cB_p(\cK)$ is a
family of
$\cB_p(\cK)$-analytic operators on $\Omega$, $\Omega\subseteq\bbC$
open. Let $\{P_n\}_{n\in\bbN}$ be a sequence of orthogonal
projections in $\cK$ converging strongly to $I_{\cK}$ as
$n\to\infty$. Then, the following limits hold uniformly with respect to
$z$ as $z$ varies in compact subsets of $\Omega$,
\begin{align}
& \lim_{n\to\infty}\|P_nL(z)P_n-L(z)\|_p=0, \lb{2.77} \\
& \lim_{n\to\infty} {\det}_p(I_{\cK}-P_nL(z)P_n)=
{\det}_p(I_{\cK}-L(z)), \lb{2.78} \\
& \lim_{n\to\infty} \f{d}{dz}{\det}_p(I_{\cK}-P_nL(z)P_n)=
\f{d}{dz}{\det}_p(I_{\cK}-L(z)). \lb{2.79}
\end{align}
\end{lemma}
%%%%%%%%%%%%%%%%%%%%%%%%%%%%%%%%%%%%%

So while the situation for analytic $\cB_p(\cK)$-valued functions
is very satisfactory, there is, however, a problem with
meromorphic (even completely meromorphic) $\cB_p(\cK)$-valued
functions as pointed out by Howland. Indeed, suppose $L(z)$,
$z\in\Omega$, is meromorphic in $\Omega$ and of finite rank. Then
of course $\det(I_{\cK}-L(\cdot))$ is meromorphic in $\Omega$.
However, the formula
\begin{equation}
{\det}_p(I_{\cK}-L(z))=\det(I_{\cK}-L(z))
\exp\bigg[\tr\bigg(-\sum_{j=1}^{p-1}j^{-1} L(z)^j\bigg)\bigg],
\quad z\in\Omega \lb{2.80}
\end{equation}
shows that ${\det}_p(I_{\cK}-L(\cdot))$, for $p>1$, in general,
will exhibit essential singularities at poles of $L$. To sidestep
this difficulty, Howland extends the definition of $m(\cdot\,;f)$ in
\eqref{2.49}, \eqref{2.49a} to functions $f$ with isolated
essential singularities as follows: Suppose $f$ is meromorphic in
$\Omega$ except at isolated essential singularities. Then we use
\eqref{2.49a} again to define
\begin{equation}
m(z;f)=\f{1}{2\pi i} \oint_{C(z; \varepsilon) } d\zeta \,
\f{f'(\zeta)}{f(\zeta)}, \quad z\in\Omega, \lb{2.81}
\end{equation}
where $\varepsilon>0$ is chosen sufficiently small to exclude all
singularities and zeros of $f$ except possibly $z$.

Given Lemma \ref{l2.14} and the extension of $m(\cdot\,;f)$ to
meromorphic functions with isolated essential singularities, Howland
\cite{Ho70} then proves the following fundamental result (the proof
of which is independent of any self-adjointness hypotheses on $H_0$
and $H$ and hence omitted here).

%%%%%%%%%%%%%%%%%%%%%%%%%%%%%%%%%%%%%
\begin{lemma} [\cite{Ho70}] \lb{l2.15}
Let $p\in\bbN$ and assume that $\{L(z)\}_{z\in\Omega}$ is a family of
$\cB_p(\cK)$-valued completely meromorphic operators on $\Omega$,
$\Omega\subseteq\bbC$ open. Let $\{M(z)\}_{z\in\Omega}$ be a boundedly
invertible operator-valued analytic function on $\Omega$ such that
\begin{equation}
M(z)[I_{\cK}-L(z)]=I_{\cK}-F(z), \quad z\in\Omega, \lb{2.82}
\end{equation}
where $F(z)$ is meromorphic and of finite rank for all $z\in\Omega$.
Define
\begin{equation}
\Delta(z)=\det(I_{\cK}-F(z)), \quad z\in\Omega, \lb{2.83}
\end{equation}
and
\begin{equation}
\Delta_p(z)={\det}_p(I_{\cK}-L(z)), \quad z\in\Omega.  \lb{2.84}
\end{equation}
Then,
\begin{equation}
m(z;\Delta)=m(z;\Delta_p), \quad z\in\Omega. \lb{2.85}
\end{equation}
\end{lemma}
%%%%%%%%%%%%%%%%%%%%%%%%%%%%%%%%%%%%%

Combining Theorem \ref{t2.11} and Lemma \ref{l2.15} yields
Howland's global Weinstein--Aronszajn formula \cite{Ho70} extended
to the non-self-adjoint case.

%%%%%%%%%%%%%%%%%%%%%%%%%%%%%%%%%%%%%
\begin{theorem} \lb{t2.16}
Assume Hypothesis \ref{h2.13}. Then the global
Weinstein--Aronszajn formula
\begin{equation}
m(z,H)=m(z,H_0)+m(z;{\det}_p(I_{\cK}-K(z))), \quad z\in\wti \rho (H_0),
\lb{2.86}
\end{equation}
holds.
\end{theorem}
%%%%%%%%%%%%%%%%%%%%%%%%%%%%%%%%%%%%%

%%%%%%%%%%%%%%%%%%%%%%%%%%%%%%%%%%%%%
\begin{remark} \lb{r2.17}
Let $H_0$ be as in Hypothesis \ref{h2.1} and assume the hypotheses on 
$\widetilde \rho(H_0)$ and $\widetilde \rho(H)$ as in Theorem \ref{t2.11} 
(cf.\ footnote~\ref{F} on p.~11). In addition, fix $p\in\bbN$, and assume
$VR_0(z)\in\cB_p(\cH)$. Define $H=H_0+V$, $\dom(H)=\dom(H_0)$. Then
\eqref{2.86} holds on $\wti \rho (H_0)$ with $K(z)=VR_0(z)$. In the
special case $p=1$ this was first obtained by Kuroda \cite{Ku61}. 
\hfill $\diamond$ 
\end{remark}
%%%%%%%%%%%%%%%%%%%%%%%%%%%%%%%%%%%%%

%%%%%%%%%%%%%%%%%%%%%%%%%%%%%%%%%%%%%%%
%%%%%%%%%%%%%%%%%%%%%%%%%%%%%%%%%%%%%%%
\section{An Application of Perturbation Determinants to \\ Schr\"odinger
Operators in Dimension $n=1,2,3$}  \lb{s6}
%%%%%%%%%%%%%%%%%%%%%%%%%%%%%%%%%%%%%%%
%%%%%%%%%%%%%%%%%%%%%%%%%%%%%%%%%%%%%%%

In dimension one on a half-line $(0,\infty)$, the perturbation determinant
associated with the Birman--Schwinger kernel corresponding to a
Schr\"odinger operator with an integrable potential on $(0,\infty)$ is
known to coincide with the corresponding Jost function and hence with a
simple Wronski determinant (cf.\ Lemmas \ref{l6.7} and \ref{l6.8}). This
reduction of an infinite-dimensional determinant to a finite-dimensional
one is quite remarkable and in this section we intend to give some ideas
as to how this fact can be generalized to dimensions two and three.

We start with the one-dimensional situation on the half-line
$\Om=(0,\infty)$ and introduce the Dirichlet and Neumann
Laplacians $H_{0,+}^D$ and $H_{0,+}^N$ in $L^2((0,\infty);dx)$ by
\begin{align}
&H_{0,+}^Df=-f'',  \no \\
&f\in \dom\big(H_{0,+}^D\big)=\{g\in L^2((0,\infty); dx) \,|\, g,g'
\in AC([0,R])
\text{ for all $R>0$}, \\
& \hspace*{6.35cm} g(0)=0, \, g''\in L^2((0,\infty); dx)\},  \no \\
&H_{0,+}^Nf=-f'',   \no \\
&f\in \dom\big(H_{0,+}^N\big)=\{g\in L^2((0,\infty); dx) \,|\, g,g'
\in AC([0,R])
\text{ for all $R>0$}, \\
& \hspace*{6.25cm} g'(0)=0, \, g''\in L^2((0,\infty); dx)\}.  \no
\end{align}

Next, we make the following assumption on the potential $V$:

%%%%%%%%%%%%%%%%%%%%%%%%%%%%%%%%%%%%%%%
\begin{hypothesis} \lb{h6.6}
Suppose $V\in L^1((0,\infty);dx)$.
\end{hypothesis}
%%%%%%%%%%%%%%%%%%%%%%%%%%%%%%%%%%%%%%%

Given Hypothesis \ref{h6.6}, we introduce the perturbed operators
$H_+^D$ and $H_+^N$ in $L^2((0,\infty);dx)$ by
\begin{align}
&H_{+}^Df=-f''+Vf,  \no \\
&f\in \dom\big(H_{0,+}^D\big)=\{g\in L^2((0,\infty); dx) \,|\, g,g'
\in AC([0,R])
\text{ for all $R>0$}, \\
& \hspace*{4.95cm} g(0)=0, \, (-g''+Vg)\in L^2((0,\infty); dx)\}, \no \\
&H_{+}^Nf=-f''+Vf,  \no \\
&f\in \dom\big(H_{0,+}^N\big)=\{g\in L^2((0,\infty); dx) \,|\, g,g'
\in AC([0,R])
\text{ for all $R>0$}, \\
& \hspace*{4.85cm} g'(0)=0, \, (-g''+Vg)\in L^2((0,\infty); dx)\}. \no
\end{align}

A fundamental system of solutions $\phi_+^D(z,\cdot)$,
$\theta_+^D(z,\cdot)$, and the Jost solution $f_+(z,\cdot)$ of
\begin{equation}
-\psi''(z,x)+V\psi(z,x)=z\psi(z,x), \quad z\in\bbC\backslash\{0\}, \;
x\geq 0,   \lb{6.5a}
\end{equation}
are introduced by
\begin{align}
\phi_+^D(z,x)&=z^{-1/2}\sin(z^{1/2}x)+\int_0^x dx' g^{(0)}_+(z,x,x')
V(x')\phi_+^D(z,x'), \\
\theta_+^D(z,x)&=\cos(z^{1/2}x)+\int_0^x dx' g^{(0)}_+(z,x,x')
V(x')\theta_+^D (z,x'), \\
f_+(z,x)&=e^{iz^{1/2}x}-\int_x^\infty dx'
g^{(0)}_+(z,x,x')V(x')f_+(z,x'),  \lb{6.8a} \\
&\hspace*{22mm} \Im(z^{1/2})\geq 0, \; z\in\bbC\backslash\{0\}, \;
x\geq 0,  \no
\end{align}
where
\begin{align}
g^{(0)}_+(z,x,x')&=z^{-1/2}\sin(z^{1/2}(x-x')).
\end{align}

We introduce
\begin{equation}
u=\exp(i\arg(V))\abs{V}^{1/2}, \quad v=\abs{V}^{1/2}, \, \text{ so
that } \, V=u\, v,
\end{equation}
and denote by $I_+$ the identity operator in $L^2((0,\infty); dx)$. In
addition, we let
\begin{equation}
W(f,g)(x)=f(x)g'(x)-f'(x)g(x), \quad x \geq 0,
\end{equation}
denote the Wronskian of $f$ and $g$, where $f,g \in C^1([0,\infty))$. We
also recall our convention to denote by $M_f$ the operator of
multiplication in $L^2((0,\infty);dx)$ by an element $f\in
L^1_{\loc}((0,\infty);dx)$ (and similarly in the higher-dimensional
context in the main part of this section).

The following is a modern formulation of a classical result by Jost and
Pais \cite{JP51}.

%%%%%%%%%%%%%%%%%%%%%%%%%%%%%%%%%%%%% 
\begin{lemma}[{\cite[Theorem\ 4.3]{GM03}}] \lb{l6.7}
Assume Hypothesis \ref{h6.6} and $z\in\bbC\backslash [0,\infty)$
with $\Im(z^{1/2})>0$. Then $\ol{M_u(H_{0,+}^D-z I_+)^{-1}M_v} \in
\cB_1(L^2((0,\infty);dx))$ and
\begin{align}
\det\big(I_+ +\ol{M_u(H_{0, +}^D-z I_+)^{-1}M_v}\big) &=
1+z^{-1/2}\int_0^\infty dx\, \sin(z^{1/2}x)V(x)f_+(z,x)   \no \\
&= W(f_+(z,\cdot),\phi_+^D(z,\cdot)) = f_+(z,0). \lb{6.75}
\end{align}
\end{lemma}
%%%%%%%%%%%%%%%%%%%%%%%%%%%%%%%%%%%%%%

Performing calculations similar to Section 4 in \cite{GM03} for the
pair of operators $H_{0,+}^N$ and $H_+^N$, one also obtains the
following result.

%%%%%%%%%%%%%%%%%%%%%%%%%%%%%%%%%%%%%
\begin{lemma} \lb{l6.8}
Assume Hypothesis \ref{h6.6} and $z\in\bbC\backslash [0,\infty)$
with $\Im(z^{1/2})>0$. Then $\ol{M_u(H_{0,+}^N-z I_+)^{-1}M_v} \in
\cB_1\big(L^2((0,\infty);dx)\big)$ and
\begin{align}
\det\big(I_+ +\ol{M_u(H_{0, +}^N-z I_+)^{-1}M_v}\big)
&= 1+ i z^{-1/2} \int_0^\infty
dx\, \cos(z^{1/2}x)V(x)f_+(z,x) \no  \\
&= - \frac{W(f_+(z,\cdot),\theta_+^D (z,\cdot))}{i z^{1/2}} =
\frac{f_+'(z,0)}{i z^{1/2}}. \lb{6.78}
\end{align}
\end{lemma}
%%%%%%%%%%%%%%%%%%%%%%%%%%%%%%%%%%%%

We emphasize that \eqref{6.75} and \eqref{6.78} exhibit the remarkable
fact that the Fredholm determinant associated with trace class operators
in the infinite-dimensional space $L^2((0,\infty); dx)$ is reduced to a
simple Wronski determinant of $\bbC$-valued distributional solutions of
\eqref{6.5a}. This fact goes back to Jost and Pais \cite{JP51} (see
also \cite{GM03}, \cite{Ne72}, \cite{Ne80}, \cite[Sect.\ 12.1.2]{Ne02},
\cite[Proposition 5.7]{Si79}, \cite{Si00}, and the extensive literature
cited in these references). The principal aim of this section is to explore
possibilities to extend this fact to higher dimensions
$n=2,3$. While a straightforward generalization of \eqref{6.75},
\eqref{6.78} appears to be difficult, we will next derive a formula for
the ratio of such determinants which permits a direct extension to
dimensions $n=2,3$.

For this purpose we introduce the boundary trace operators
$\ga_D$ (Dirichlet trace) and $\ga_N$ (Neumann trace) which, in the
current one-dimensional half-line situation, are just the functionals,
\begin{equation}
\ga_D \colon \begin{cases} C([0,\infty)) \to \bbC \\
\hspace*{1.3cm} g \mapsto g(0) \end{cases}, \quad
\ga_N \colon \begin{cases}C^1([0,\infty)) \to \bbC  \\
\hspace*{1.43cm} h \mapsto   - h'(0) \end{cases}.
\end{equation}
In addition, we denote by $m_{0,+}^D$, $m_+^D$, $m_{0,+}^N$, and $m_+^N$
the Weyl--Titchmarsh $m$-functions corresponding to $H_{0,+}^D$,
$H_{+}^D$, $H_{0,+}^N$, and $H_{+}^N$, respectively,
\begin{align}
&m_{0,+}^D(z) = i z^{1/2},\ \quad m_{0,+}^N (z)= -\frac{1}{m_{0,+}^D(z)}
= i z^{-1/2},  \lb{6.16a} \\
&m_{+}^D(z) = \frac{f_+'(z,0)}{f_+(z,0)}, \quad m_{+}^N (z)=
-\frac{1}{m_{+}^D(z)} = -\frac{f_+(z,0)}{f_+'(z,0)}.  \lb{6.17a}
\end{align}

%%%%%%%%%%%%%%%%%%%%%%%%%%%%%%%%%%%%%
\begin{theorem} \lb{t6.9}
Assume Hypothesis \ref{h6.6} and let $z\in\bbC\backslash\si(H_+^D)$
with $\Im(z^{1/2})>0$. Then,
\begin{align}
& \frac{\det\big(I_+ +\ol{M_u(H_{0, +}^N-z I_+)^{-1} M_v}\big)}
{\det\big(I_+ +\ol{M_u(H_{0, +}^D-z I_+)^{-1}M_v}\big)}
\no \\
& \quad = \f{W(f_+(z),\phi_+^N(z))}{i z^{1/2}W(f_+(z),\phi_+^D(z))} =
\f{f'_+(z,0)}{i z^{1/2}f_+(z,0)} = \f{m_+^D(z)}{m_{0,+}^D(z)} =
\f{m_{0,+}^N(z)}{m_+^N(z)} \lb{6.82}
\\
&\quad = 1 - \big(\ol{\ga_N(H_+^D-z I_+)^{-1}M_V
\big[\ga_D(H_{0,+}^N-\ol{z}I_+)^{-1}\big]^*}\big) 1. \lb{6.83}
\end{align}
\end{theorem}
%%%%%%%%%%%%%%%%%%%%%%%%%%%%%%%%%%%%
\begin{proof}
We start by noting that $\si(H_{0,+}^D)=\si(H_{0,+}^N)=[0,\infty)$.
Applying Lemmas \ref{l6.7} and \ref{l6.8} and equations \eqref{6.16a}
and \eqref{6.17a} proves \eqref{6.82}.

To verify the equality of \eqref{6.82} and \eqref{6.83} requires some
preparations. First we recall that the Green's functions (i.e., integral
kernels) of the resolvents of $H^D_{0,+}$ and $H^N_{0,+}$ are given by
\begin{align}
(H^D_{0,+} -zI_+)^{-1}(x,x')&=\begin{cases} \f{\sin(z^{1/2}x)}{z^{1/2}}
e^{i z^{1/2}x'}, & 0\leq x\leq x', \\
\f{\sin(z^{1/2}x')}{z^{1/2}}
e^{i z^{1/2}x}, & 0\leq x' \leq x, \end{cases}   \\
(H^N_{0,+} -zI_+)^{-1}(x,x')&=\begin{cases} \f{\cos(z^{1/2}x)}{-i z^{1/2}}
e^{i z^{1/2}x'}, & 0\leq x\leq x', \\
\f{\cos(z^{1/2}x')}{-i z^{1/2}}
e^{i z^{1/2}x}, & 0\leq x' \leq x, \end{cases}
\end{align}
and hence Krein's formula for the resolvent difference of $H^D_{0,+}$ and
$H^N_{0,+}$ takes on the simple form
\begin{align}
&(H^D_{0,+} -zI_+)^{-1}-(H^N_{0,+} -zI_+)^{-1}=-i z^{-1/2}
(\ol{\psi_{0,+}(z,\cdot)}, \cdot)_{L^2((0,\infty);dx)}\psi_{0,+}(z,\cdot),
\no  \\
& \hspace*{5.2cm}  z\in\rho(H^D_{0,+})\cap\rho(H^N_{0,+}), \;
\Im(z^{1/2})>0,   \lb{6.22a}
\end{align}
where we abbreviated
\begin{equation}
\psi_{0,+}(z,x)= e^{i z^{1/2}x}, \quad \Im(z^{1/2})>0, \; x\geq 0.
\end{equation}
We also recall
\begin{equation}
(H^D_{+} -zI_+)^{-1}(x,x')=\begin{cases} \phi^D_+(z,x) \psi_+(z,x'), &
0\leq x\leq x', \\
\phi^D_+(z,x') \psi_+(z,x), & 0\leq x' \leq x, \end{cases}
\end{equation}
where
\begin{equation}
\psi_+(z,x)=\theta^D_+(z,x) + m^D_+(z) \phi^D_+(z,x), \quad
z\in\rho(H^D_+), \; x\geq 0,
\end{equation}
and
\begin{equation}
\psi_+(z,\cdot)=\f{f_+(z,\cdot)}{f_+(z,0)} \in L^2((0,\infty);dx), \quad
z\in\rho(H^D_+).
\end{equation}
In fact, a standard iteration argument applied to \eqref{6.8a} shows that
\begin{equation}
|\psi_+(z,x)|\leq C(z) e^{-\Im(z^{1/2})x}, \quad 
\Im(z^{1/2})> 0, \; x\geq 0.
\end{equation}
In addition, we note that
\begin{align}
\gamma_N(H^D_{0,+} -zI_+)^{-1} g &= -\int_0^\infty dx \,
e^{i z^{1/2}x} g(x), \quad g\in L^2((0,\infty);dx),
\\
\ga_N(H_+^D-zI_+)^{-1}g &= -\int_0^\infty dx\, \psi_+(z,x)g(x), \quad
g\in L^2((0,\infty);dx),
\\
\gamma_D(H^N_{0,+} -zI_+)^{-1} f &= i z^{-1/2}\int_0^\infty dx\,
e^{i z^{1/2}x} f(x), \quad f\in L^2((0,\infty);dx),
\end{align}
and hence,
\begin{equation}
\big(\big[\ga_D(H^N_{0,+}-\ol{z}I_+)^{-1}\big]^* \,
c\big)(\cdot)=ic z^{-1/2}\psi_{0,+}(z,\cdot), \quad c\in\bbC.
\end{equation}
Then Krein's formula \eqref{6.22a}  can be rewritten as
\begin{align}
&(H^D_{0,+} -zI_+)^{-1}-(H^N_{0,+} -zI_+)^{-1}=
\big[\ga_D(H^N_{0,+}-\ol{z}I_+)^{-1}\big]^* \gamma_N(H^D_{0,+}
-zI_+)^{-1},
\no  \\
& \hspace*{5.2cm}  z\in\rho(H^D_{0,+})\cap\rho(H^N_{0,+}), \;
\Im(z^{1/2})>0.   \lb{6.27a}
\end{align}
Finally, using the facts (cf.\ \eqref{6.8a})
\begin{align}
f_+(z,0)&=1+ z^{-1/2}
\int_0^\infty dx \, \sin(z^{1/2}x) V(x) f_+(z,x),  \\
f'_+(z,0)&=i z^{1/2}- \int_0^\infty dx\, \cos(z^{1/2}) V(x) f_+(z,x),
\end{align}
one computes (since $v\in L^2(\bbR;dx)$ and $\psi_+(z,\cdot)\in
L^\infty(\bbR;dx)$)
\begin{align}
&-\big[\ol{\gamma_N (H^D_{+} -zI_+)^{-1} M_V
\big[\ga_D(H^N_{0,+}-\ol{z}I_+)^{-1}\big]^*}\big] 1
\no \\
& \quad = -i z^{-1/2}\ol{\gamma_N (H^D_{+} -zI_+)^{-1} M_u} (v
\psi_{0,+})(z,\cdot)  \no \\
& \quad = i z^{-1/2} \int_0^\infty dx \, e^{i z^{1/2}x} V(x)
\psi_+(z,x)  \no \\
& \quad = i z^{-1/2} \int_0^\infty dx \, \bigg[\cos(z^{1/2}x)
+i z^{1/2} \f{\sin(z^{1/2}x)}{z^{1/2}}\bigg] V(x)
\f{f_+(z,x)}{f_+(z,0)}  \no \\
& \quad =
\f{i}{z^{1/2} f_+(z,0)}[i z^{1/2} - f'_+(z,0) +i z^{1/2}(f_+(z,0)-1)]
\no \\
& \quad = \f{f'_+(z,0)}{i z^{1/2} f_+(z,0)} -1.
\end{align}
\end{proof}
%%%%%%%%%%%%%%%%%%%%%%%%%%%%%%%%%%%%%%

At first sight it may seem unusual to even attempt to prove \eqref{6.83}
in the one-dimensional case since \eqref{6.82} already yields the
reduction of a Fredholm determinant to a simple Wronski determinant.
However, we will see in Theorem \ref{t6.5} that it is precisely
\eqref{6.83} that permits a straightforward extension to dimensions
$n=2,3$.

%%%%%%%%%%%%%%%%%%%%%%%%%%%%%%%%%%%%%
\begin{remark} \lb{r6.5}
As in Theorem \ref{t6.9} we assume Hypothesis \ref{h6.6} and suppose
$z\in\bbC\backslash\si(H_+^D)$. First we note
that
\begin{align}
&(H_{0,+}^D-zI_+)^{-1/2}(H_+^D-zI_+)(H_{0,+}^D-zI_+)^{-1/2} - I_+
\in\cB_1\big(L^2((0,\infty);dx)\big), \lb{6.79} \\
&(H_{0,+}^N-zI_+)^{-1/2}(H_+^N-zI_+)(H_{0,+}^N-zI_+)^{-1/2} - I_+
\in\cB_1\big(L^2((0,\infty);dx)\big). \lb{6.80}
\end{align}
Indeed, it follows from the proof of
\cite[Theorem 4.2]{GM03} (cf.\ also Lemma \ref{l6.2} below), that
\begin{align}
& \ol{(H_{0,+}^D-z I_+)^{-1/2}M_u}, \; M_v(H_{0,+}^D-z I_+)^{-1/2} \in
\cB_2\big(L^2((0,\infty);dx)\big),
\end{align}
and hence,
\begin{align}
&(H_{0,+}^D-zI_+)^{-1/2}(H_+^D-z I_+)(H_{0,+}^D-z I_+)^{-1/2} - I_+
\\
&\quad = (H_{0,+}^D-zI_+)^{-1/2}M_V(H_{0,+}^D-zI_+)^{-1/2}
\in\cB_1\big(L^2((0,\infty);dx)\big).
\end{align}
This proves \eqref{6.79}, and a similar argument yields \eqref{6.80}.
Using the cyclicity of $\det(\cdot)$, one can then rewrite the left-hand
side of \eqref{6.82} as follows,
\begin{align}
& \frac{\det\big(I_+ +\ol{M_u(H_{0, +}^N-z I_+)^{-1} M_v}\big)}
{\det\big(I_+ +\ol{M_u(H_{0, +}^D-z I_+)^{-1}M_v}\big)} \no \\
&\quad = \frac{\det\big(
I_++(H_{0,+}^N-zI_+)^{-1/2}M_V(H_{0,+}^N-zI_+)^{-1/2} \big)}
{\det\big(I_++(H_{0,+}^D-zI_+)^{-1/2}M_V(H_{0,+}^D-zI_+)^{-1/2}\big)}
\no
\\
& \quad =\frac{\det\big(
(H_{0,+}^N-zI_+)^{-1/2}(H_+^N-zI_+)(H_{0,+}^N-zI_+)^{-1/2} \big)}
{\det\big( (H_{0,+}^D-zI_+)^{-1/2}(H_+^D-zI_+)(H_{0,+}^D-zI_+)^{-1/2}
\big)}.  \lb{6.87}
\end{align}
Equation \eqref{6.87} illustrates the kind of symmetrized
perturbation determinants underlying Theorem \ref{t6.9}. 
\hfill $\diamond$
\end{remark}
%%%%%%%%%%%%%%%%%%%%%%%%%%%%%%%%%%%%

Now we turn to dimensions $n=2,3$. As a general rule, we will have
to replace Fredholm determinants by modified ones.

For the remainder of this section we make the following assumptions on the
domain $\Omega\subset \bbR^n$, $n=2,3$, and the potential $V$:

%%%%%%%%%%%%%%%%%%%%%%%%%%%%%%%%%%%%
\begin{hypothesis} \lb{h6.1} Let $n=2,3$. \\
$(i)$ Assume that $\Omega\subset{\bbR}^n$ is an open nonempty domain of
class $C^{1,r}$ for some $(1/2)<r<1$ with a compact, nonempty boundary,
$\partial\Omega$. (For details we refer to Appendix \ref{sA}.) \\
$(ii)$ Suppose that $V\in L^2(\Om;d^nx)$. \\
\end{hypothesis}
%%%%%%%%%%%%%%%%%%%%%%%%%%%%%%%%%%%%

First we introduce the boundary trace operator $\ga_D^0$  (Dirichlet
trace) by
\begin{equation}
\ga_D^0\colon C(\ol{\Om})\to C(\dOm), \quad
\ga_D^0 u = u|_\dOm .
\end{equation}
Then there exists a bounded, linear operator $\gamma_D$,
\begin{equation}
\ga_D\colon H^{s}(\Om)\to H^{s-1/2}(\dOm) \hookrightarrow
L^2(\dOm;d^{n-1}\si), \quad 1/2<s<3/2,
\lb{6.1}
\end{equation}
whose action is compatible with $\ga_D^0$, that is, the two Dirichlet
trace  operators coincide on the intersection of their domains.
It is well-known (see, e.g., \cite[Theorem 3.38]{Mc00}), that $\ga_D$ is
bounded. Here $d^{n-1}\sigma$ denotes the surface measure
on $\dOm$ and we refer to Appendix \ref{sA} for our notation in connection
with Sobolev spaces.

Next, let $I_{\partial\Om}$ denote the identity operator in
$L^2(\partial\Om;d^{n-1}\sigma)$, and introduce the operator $\ga_N$
(Neumann trace) by
\begin{align}
\ga_N = \nu\cdot\ga_D\nabla \colon H^{s+1}(\Om)\to L^2(\dOm;d^{n-1}\si),
\quad 1/2<s<3/2, \lb{6.2}
\end{align}
where $\nu$ denotes outward pointing normal unit vector to $\partial\Om$.
It follows from \eqref{6.1} that $\ga_N$ is also a bounded operator.

Given Hypothesis \ref{h6.1}\,$(i)$, we introduce the Dirichlet and Neumann
Laplacians $H_{0,\Om}^D$ and $H_{0,\Om}^N$ associated with the domain
$\Om$ as follows,
\begin{align}
H^D_{0,\Om} = -\Delta, \quad \dom(H^D_{0,\Om}) = \{u\in H^{2}(\Om)
\,|\, \ga_D u = 0\}, \lb{6.3}
\\
H^N_{0,\Om} = -\Delta, \quad \dom(H^N_{0,\Om}) = \{u\in H^{2}(\Om)
\,|\, \ga_N u = 0\}. \lb{6.4}
\end{align}

In the following we denote by $I_\Om$ the identity operator
in $L^2(\Om;d^nx)$.

%%%%%%%%%%%%%%%%%%%%%%%%%%%%%%%%%%%%%
\begin{lemma} \lb{l6.1}
Assume Hypothesis \ref{h6.1}\,$(i)$. Then the operators $H_{0,\Om}^D$ and
$H_{0,\Om}^N$ introduced in \eqref{6.3} and \eqref{6.4} are nonnegative
and self-adjoint in $\cH=L^2(\Om;d^nx)$ and the following mapping
properties hold for all %$q\in\bbR$ 
$q\in [0,1]$ and $z\in\bbC\backslash[0,\infty)$,
\begin{align}
(H_{0,\Om}^D-zI_{\Om})^{-q},\, (H_{0,\Om}^N-zI_{\Om})^{-q}
\in\cB\big(L^2(\Om;d^nx),H^{2q}(\Om)\big). \lb{6.5}
\end{align}
\end{lemma}
%%%%%%%%%%%%%%%%%%%%%%%%%%%%%%%%%%%%

The fractional powers in \eqref{6.5} (and in subsequent analogous cases
such as in \eqref{6.51a}) are defined via the functional calculus implied
by the spectral theorem for self-adjoint operators. For the proof of Lemma
\ref{l6.1} we refer to Lemmas \ref{lA.1} and
\ref{lA.2} in Appendix \ref{sA}.

%%%%%%%%%%%%%%%%%%%%%%%%%%%%%%%%%%%%
\begin{lemma} \lb{l6.2}
Assume Hypothesis \ref{h6.1}\,$(i)$ and let $(n/2p)<q\leq1$, $p\geq
2$, $n=2, 3$, $f\in L^p(\Om;d^nx)$, and
$z\in\bbC\backslash[0,\infty)$. Then,
\begin{align} \lb{6.11}
M_f(H_{0,\Om}^D-zI_{\Om})^{-q}, \, M_f(H_{0,\Om}^N-zI_{\Om})^{-q}
\in\cB_p\big(L^2(\Om;d^nx)\big)
\end{align}
and for some $c(z)>0$ $($independent of $f$\,$)$
\begin{align}
\begin{split}
&\norm{M_f(H_{0,\Om}^D-zI_{\Om})^{-q}}_{\cB_p(L^2(\Om;d^nx))}+
\norm{M_f(H_{0,\Om}^N-zI_{\Om})^{-q}}_{\cB_p(L^2(\Om;d^nx))}   \\
& \quad \leq c(z)\, \|(\abs{\cdot}^2-z)^{-q}\|_{L^p(\bbR^n;d^nx)}
\|f\|_{L^p(\Om;d^nx)}. \lb{6.12}
\end{split}
\end{align}
\end{lemma}
%%%%%%%%%%%%%%%%%%%%%%%%%%%%%%%%%%%
\begin{proof}
We start by noting that under the assumption that $\Om$ is a Lipschitz
domain, there is a bounded extension operator $\cE$,
\begin{equation}
\cE\in\cB\big(H^{2q}(\Om),H^{2q}(\bbR^n)\big) \,
\text{ such that } \, (\cE u){|_\Om} = u, \quad u\in H^{2q}(\Om)
\lb{6.43a}
\end{equation}
(see, e.g., \cite[Theorem\ A.4]{Mc00}). Next, denote by $\cR_{\Om}$ the
restriction operator
\begin{equation}
\cR_{\Om}\colon \begin{cases} L^2(\bbR^n;d^nx) \to L^2(\Om;d^nx), \\
\hspace*{1.65cm}u \mapsto u|_{\Om}, \end{cases}
\end{equation}
and let $\ti f$ denote the following extension of $f$,
\begin{equation}
\ti f(x) =
\begin{cases}f(x), & x\in\Om, \\ 0, &
x\in\bbR^n\backslash\Om, \end{cases} \quad  \ti f\in
L^p(\bbR^n;d^nx).
\end{equation}
Then,
\begin{equation} \lb{6.14}
M_f (H_{0,\Om}^D-zI_{\Om})^{-q}= \cR_\Om M_{\ti f}
(H_{0}-zI)^{-q}(H_{0}-zI)^{q}\cE (H_{0,\Om}^D-zI_{\Om})^{-q},
\end{equation}
where (for simplicity) $I$ denotes the identity operator in
$L^2(\bbR^n;d^nx)$ and $H_0$ denotes the nonnegative self-adjoint operator
\begin{equation}
H_0 = -\Delta, \quad \dom(H_0)=H^{2}(\bbR^n)
\end{equation}
in $L^2(\bbR^n;d^nx)$. Utilizing the representation of $(H_{0}-zI)^{q}$ as
the operator of multiplication by $\big(\abs{\xi}^2-z\big)^{q}$ in the
Fourier space $L^2(\bbR^n;d^n\xi)$, one obtains
\begin{equation}
(H_{0}-zI)^{q}\in\cB\big(H^{2q}(\bbR^n),L^2(\bbR^n;d^nx)\big),  \lb{6.51a}
\end{equation}
which together with \eqref{6.5} and the mapping property of the
extension operator $\cE$ in \eqref{6.43a} yields
\begin{equation}
(H_{0}-zI)^{q}\cE (H_{0,\Om}^D-zI_{\Om})^{-q}\in
\cB\big(L^2(\Om;d^nx),L^2(\bbR^n;d^nx)\big). \lb{6.15}
\end{equation}
By  \cite[Theorem 4.1]{Si79} (or \cite[Theorem XI.20]{RS79}) one also
obtains
\begin{align}
M_{\ti f}(H_{0}-zI)^{-q}\in\cB_p\big(L^2(\bbR^n;d^nx)\big) \lb{6.18}
\end{align}
and
\begin{align}\begin{split}
\big\|M_{\ti f}(H_{0}-zI)^{-q}\big\|_{\cB_p(L^2(\bbR^n;d^nx))} &\leq c\,
\|(\abs{\cdot}^2-z)^{-q}\|_{L^p(\bbR^n;d^nx)} \|\ti
f\|_{L^p(\bbR^n;d^nx)} \\&= c\,
\|(\abs{\cdot}^2-z)^{-q}\|_{L^p(\bbR^n;d^nx)}
\|f\|_{L^p(\Om;d^nx)}. \lb{6.19}
\end{split}\end{align}
Thus, the Dirichlet parts of \eqref{6.11} and \eqref{6.12} follow from
\eqref{6.14}, \eqref{6.15}, \eqref{6.18}, and \eqref{6.19}.

Similar arguments prove the Neumann parts of \eqref{6.11} and
\eqref{6.12}.
\end{proof}
%%%%%%%%%%%%%%%%%%%%%%%%%%%%%%%%%%%

%%%%%%%%%%%%%%%%%%%%%%%%%%%%%%%%%%%
\begin{lemma} \lb{l6.3}
Assume Hypothesis \ref{h6.1}\,$(i)$ and let $\eps\in(0,1]$, $n=2,3$, and
$z\in\bbC\backslash[0,\infty)$. Then,
\begin{align}
\ga_N(H_{0,\Om}^D-zI_{\Om})^{-\frac{3+\eps}{4}},
\ga_D(H_{0,\Om}^N-zI_{\Om})^{-\frac{1+\eps}{4}} \in
\cB\big(L^2(\Om;d^nx),L^2(\dOm;d^{n-1}\si)\big).  \lb{6.47a}
\end{align}
\end{lemma}
%%%%%%%%%%%%%%%%%%%%%%%%%%%%%%%%%%%%
\begin{proof}
It follows from \eqref{6.5}, that
\begin{align}
(H_{0,\Om}^D-zI_{\Om})^{-\frac{3+\eps}{4}} &\in
\cB\big(L^2(\Om;d^nx),H^{\f{3+\eps}{2}}(\Om)\big),
\\
(H_{0,\Om}^N-zI_{\Om})^{-\frac{1+\eps}{4}} &\in
\cB\big(L^2(\Om;d^nx),H^{\f{1+\eps}{2}}(\Om)\big),
\end{align}
and hence one infers the result from \eqref{6.1} and \eqref{6.2}.
\end{proof}
%%%%%%%%%%%%%%%%%%%%%%%%%%%%%%%%%%%%

%%%%%%%%%%%%%%%%%%%%%%%%%%%%%%%%%%%%
\begin{corollary} \lb{c6.4}
Let $f_1\in L^{p_1}(\Om;d^nx)$, $p_1>2n$, $f_2\in L^{p_2}(\Om;d^nx)$,
$p_2\geq 2$, $p_2>2n/3$, $n=2, 3$, and $z\in\bbC\backslash[0,\infty)$.
Then,
\begin{align}
\ol{\ga_N(H_{0,\Om}^D-zI_{\Om})^{-1}M_{f_1}} &\in
\cB_{p_1}\big(L^2(\Om;d^nx),L^2(\dOm;d^{n-1}\si)\big),   \lb{6.50}
\\
\ol{\ga_D(H_{0,\Om}^N-zI_{\Om})^{-1}M_{f_2}} &\in
\cB_{p_2}\big(L^2(\Om;d^nx),L^2(\dOm;d^{n-1}\si)\big)   \lb{6.51}
\end{align}
and for some $c_j(z)>0$ $($independent of $f_j$$)$, $j=1,2$,
\begin{align}
\norm{\ol{\ga_N(H_{0,\Om}^D-zI_{\Om})^{-1}M_{f_1}}}_
{\cB_{p_1}(L^2(\Om;d^nx),L^2(\dOm;d^{n-1}\si))} &\leq
c_1(z) \norm{f_1}_{L^{p_1}(\Om;d^nx)}, \lb{6.25}
\\
\norm{\ol{\ga_D(H_{0,\Om}^N-zI_{\Om})^{-1}M_{f_2}}}_
{\cB_{p_2}(L^2(\Om;d^nx),L^2(\dOm;d^{n-1}\si))} &\leq
c_2(z) \norm{f_2}_{L^{p_2}(\Om;d^nx)}. \lb{6.26}
\end{align}
\end{corollary}
%%%%%%%%%%%%%%%%%%%%%%%%%%%%%%%%%%% 
\begin{proof}
Let $\eps_1,\eps_2\in (0,1)$ be such that $0<\eps_1<1-(2n/p_1)$ and
$0<\eps_2<\min\{1,3-(2n/p_2)\}$. Then,
\begin{align}
\ol{\ga_N(H_{0,\Om}^D-zI_{\Om})^{-1}M_{f_1}} &=
\ga_N(H_{0,\Om}^D-zI_{\Om})^{-\frac{3+\eps_1}{4}}
\ol{(H_{0,\Om}^D-zI_{\Om})^{-\frac{1-\eps_1}{4}}M_{f_1}},
\\
\ol{\ga_D(H_{0,\Om}^N-zI_{\Om})^{-1}M_{f_2}} &=
\ga_D(H_{0,\Om}^N-zI_{\Om})^{-\frac{1+\eps_2}{4}}
\ol{(H_{0,\Om}^N-zI_{\Om})^{-\frac{3-\eps_2}{4}}M_{f_2}},
\end{align}
together with Lemmas \ref{l6.2} and \ref{l6.3} prove the corollary.
\end{proof}
%%%%%%%%%%%%%%%%%%%%%%%%%%%%%%%%%%%%%

Next, we introduce the perturbed operators $H_{\Om}^D$ and
$H_{\Om}^N$ in $L^{2}(\Om;d^nx)$ as follows. We denote by $A=M_u$ and
$B=B^*=M_v$ the operators of multiplication by
$u=\exp(i\arg(V))\abs{V}^{1/2}$ and $v=\abs{V}^{1/2}$ in
$L^{2}(\Om;d^nx)$, respectively, so that $M_V=B A=M_u M_v$. Applying
Lemma \ref{l6.2} to $f=u\in L^4(\Om;d^nx)$ with $q=1/2$ yields
\begin{align}
M_u(H_{0,\Om}^D-zI_{\Om})^{-1/2}, \, \ol{(H_{0,\Om}^D-zI_{\Om})^{-1/2}M_v}
&\in\cB_4\big(L^2(\Om;d^nx)\big), \quad z\in\bbC\backslash [0,\infty),
\lb{6.6}
\\
M_u(H_{0,\Om}^N-zI_{\Om})^{-1/2}, \, \ol{(H_{0,\Om}^N-zI_{\Om})^{-1/2}M_v}
&\in\cB_4\big(L^2(\Om;d^nx)\big), \quad z\in\bbC\backslash [0,\infty),
\lb{6.7}
\end{align}
and hence in particular,
\begin{align}
& \dom(A)=\dom(B) \supseteq  H^{1}(\Om) \supset H^{2}(\Om)
\supseteq \dom(H^N_{0,\Om}),  \lb{6.56a} \\
& \dom(A)=\dom(B) \supseteq  H^{1}(\Om) \supseteq H^{1}_0(\Om)
\supseteq \dom(H^D_{0,\Om}).  \lb{6.57a}
\end{align}
Thus, Hypothesis \ref{2.1}\,$(i)$ is satisfied for $H_{0,\Om}^D$ and
$H_{0,\Om}^N$. Moreover, \eqref{6.6} and \eqref{6.7} imply
\begin{align}
&\ol{M_u(H_{0,\Om}^D-zI_{\Om})^{-1}M_v}, \,
\ol{M_u(H_{0,\Om}^N-zI_{\Om})^{-1}M_v}
\in\cB_2\big(L^2(\Om;d^nx)\big), \quad z\in\bbC\backslash [0,\infty),
\lb{6.8}
\end{align}
which verifies Hypothesis \ref{2.1}\,$(ii)$ for $H_{0,\Om}^D$ and
$H_{0,\Om}^N$. One verifies Hypothesis \ref{2.1}\,$(iii)$ by
utilizing \eqref{6.12} with 
$-z>0$ sufficiently large, such that the $\cB_4$-norms of the operators in \eqref{6.6}
and \eqref{6.7} are less than 1, and hence, the Hilbert--Schmidt
norms of the operators in \eqref{6.8} are less than 1. Thus, applying
Theorem \ref{2.3} one obtains the densely defined, closed
operators $H_{\Om}^D$ and $H_{\Om}^N$ (which are extensions of
$H_{0,\Om}^D+M_V$ on $\dom(H_{0,\Om}^D)\cap\dom(M_V)$ and
$H_{0,\Om}^N+M_V$ on
$\dom(H_{0,\Om}^N)\cap\dom(M_V)$, respectively).

We note in passing that \eqref{6.5}--\eqref{6.12},
\eqref{6.47a}, \eqref{6.50}--\eqref{6.26}, \eqref{6.6}--\eqref{6.8}, etc.,
extend of course to all $z$ in the resolvent set of the corresponding
operators $H_{0,\Om}^D$ and $H_{0,\Om}^N$.

The following result is a direct extension of the one-dimensional result
in Theorem \ref{t6.9}.

%%%%%%%%%%%%%%%%%%%%%%%%%%%%%%%%%%%%%
\begin{theorem} \lb{t6.5}
Assume Hypothesis \ref{h6.1} and
$z\in\bbC\backslash\big(\si(H_{\Om}^D)\cup
\si(H_{0,\Om}^D) \cup \si(H_{0,\Om}^N)\big)$. Then,
\begin{align}
&\ol{\ga_N(H_{0,\Om}^D-zI_{\Om})^{-1}M_V
(H_{\Om}^D-zI_{\Om})^{-1}M_V
\big[\ga_D(H_{0,\Om}^N-\ol{z}I_{\Om})^{-1}\big]^*}
\no  \\
& \quad
\in\cB_1\big(L^2(\dOm;d^{n-1}\si)\big),  \lb{6.70a}  \\
&\ol{\ga_N(H_{\Om}^D-zI_{\Om})^{-1}M_V
\big[\ga_D(H_{0,\Om}^N-\ol{z}I_{\Om})^{-1}\big]^*}
\in\cB_2\big(L^2(\dOm;d^{n-1}\si)\big), \lb{6.70b}
\end{align}
and
\begin{align}
& \frac{\det{}_2\big(I_{\Om}+\ol{M_u(H_{0,\Om}^N-zI_{\Om})^{-1}M_v}\big)}
{\det{}_2\big(I_{\Om}+\ol{M_u(H_{0,\Om}^D-zI_{\Om})^{-1}M_v}\big)}  \no \\
&\quad = \det{}_2\big(I_{\dOm} -
\ol{\ga_N(H_{\Om}^D-zI_{\Om})^{-1}M_V
\big[\ga_D(H_{0,\Om}^N-\ol{z}I_{\Om})^{-1}\big]^*}\big)
\lb{6.31} \\
&\quad\quad \times \exp\big(\tr\big(
\ol{\ga_N(H_{0,\Om}^D-zI_{\Om})^{-1}M_V
(H_{\Om}^D-zI_{\Om})^{-1}M_V
\big[\ga_D(H_{0,\Om}^N-\ol{z}I_{\Om})^{-1}\big]^*}
\big)\big). \no
\end{align}
\end{theorem}
%%%%%%%%%%%%%%%%%%%%%%%%%%%%%%%%%%%%%% 
\begin{proof}
From the outset we note that the left-hand side of \eqref{6.31} is
well-defined by \eqref{6.8}. Let $z\in\bbC\backslash\big(\si(H_{\Om}^D)\cup
\si(H_{0,\Om}^D) \cup \si(H_{0,\Om}^N)\big)$ and
\begin{align}
u(x) &= \exp(i\arg(V(x)))\abs{V(x)}^{1/2},\quad
v(x)=\abs{V(x)}^{1/2},
\\
\wti u(x) &= \exp(i\arg(V(x)))\abs{V(x)}^{5/6},\quad \wti
v(x)=\abs{V(x)}^{1/6}.
\end{align}
Next, we introduce
\begin{equation}
K_D(z)=-\ol{M_u(H_{0,\Om}^D-zI_{\Om})^{-1}M_v}, \quad
K_N(z)=-\ol{M_u(H_{0,\Om}^N-zI_{\Om})^{-1}M_v}
\end{equation}
(cf.\ \eqref{2.4}) and utilize the following facts,
\begin{align}
[I_{\Om}-K_D(z)]^{-1} &= I_{\Om} + K_D(z)[I_{\Om}-K_D(z)]^{-1},
\\
[I_{\Om}-K_D(z)]^{-1} &\in\cB(L^2(\Om;d^nx)),
\end{align}
and
\begin{align}
1& = \det{}_2(I_{\Om}) = \det{}_2\big([I_{\Om}-K_D(z)][I-K_D(z)]^{-1}\big)
\\
&= \det{}_2\big(I_{\Om}-K_D(z)\big)\det{}_2\big([I_{\Om}-K_D(z)]^{-1}\big)
\exp\big(\tr\big(K_D(z)^2[I_{\Om}-K_D(z)]^{-1}\big)\big).  \no
\end{align}
Thus, one obtains
\begin{align}
&\det{}_2\big([I_{\Om}-K_N(z)][I_{\Om}-K_D(z)]^{-1}\big)  \no \\
&\quad = \det{}_2\big(I_{\Om}-K_N(z)\big)
\det{}_2\big([I_{\Om}-K_D(z)]^{-1}\big) \no \\
& \qquad \, \times
\exp\big(\tr\big(K_N(z)K_D(z)[I_{\Om}-K_D(z)]^{-1}\big)\big) \\
&\quad =
\frac{\det{}_2\big(I_{\Om}-K_N(z)\big)}{\det{}_2\big(I_{\Om}-K_D(z)\big)}
\exp\big(\tr\big((K_N(z)-K_D(z))K_D(z)[I_{\Om}-K_D(z)]^{-1}\big)\big). \no
\end{align}
At this point, the left-hand side of \eqref{6.31} can be rewritten as
\begin{align}
&\frac{\det{}_2\big(I_{\Om}+\ol{M_u(H_{0,\Om}^N-zI_{\Om})^{-1}M_v}\big)}
{\det{}_2\big(I_{\Om}+\ol{M_u(H_{0,\Om}^D-zI_{\Om})^{-1}M_v}\big)} =
\frac{\det{}_2\big(I_{\Om}-K_N(z)\big)}{\det{}_2\big(I_{\Om}-K_D(z)\big)}
\no \\
&\quad = \det{}_2\big([I_{\Om}-K_N(z)][I_{\Om}-K_D(z)]^{-1}\big)
\no \\
&\quad\quad \, \times
\exp\big(\tr\big((K_D(z)-K_N(z))K_D(z)[I_{\Om}-K_D(z)]^{-1}\big)\big)
\no  \\
&\quad =
\det{}_2\big(I_{\Om}+(K_D(z)-K_N(z))[I_{\Om}-K_D(z)]^{-1}\big)
\lb{6.43}
\\
&\quad\quad \, \times
\exp\big(\tr\big((K_D(z)-K_N(z))K_D(z)[I_{\Om}-K_D(z)]^{-1}\big)\big).
\no
\end{align}
Next, temporarily suppose that $V\in L^2(\Om;d^nx)\cap
L^6(\Om;d^nx)$. Using Lemma \ref{lA.3} (an extension of a result of
Nakamura \cite[Lemma 6]{Na01}) and Remark \ref{rA.4}, one finds
\begin{align}
\begin{split}
K_D(z)-K_N(z) &=
-\ol{M_u\big[(H_{0,\Om}^D-zI_{\Om})^{-1}-(H_{0,\Om}^N
-zI_{\Om})^{-1}\big]M_v} \lb{6.44}
\\ &=
-\ol{M_u\big[\ga_D(H_{0,\Om}^N-\ol{z}I_{\Om})^{-1}\big]^*}\,
\ol{\ga_N(H_{0,\Om}^D -zI_{\Om})^{-1}M_v},
\\ &=
-\big[\ol{\ga_D(H_{0,\Om}^N-\ol{z}I_{\Om})^{-1}M_{\ol{u}}}\big]^*
\ol{\ga_N(H_{0,\Om}^D -zI_{\Om})^{-1}M_v}.
\end{split}
\end{align}
Thus, inserting \eqref{6.44} into \eqref{6.43} yields,
\begin{align}
&\frac{\det{}_2\big(I_{\Om}+\ol{M_u(H_{0,\Om}^N-zI_{\Om})^{-1}M_v}\big)}
{\det{}_2\big(I_{\Om}+\ol{M_u(H_{0,\Om}^D-zI_{\Om})^{-1}M_v}\big)}  \no
\\
&\quad = \det{}_2\Big(I_{\Om} - \big[\ol{\ga_D
(H_{0,\Om}^N-\ol{z}I_{\Om})^{-1}M_{\ol{u}}}\big]^*
\ol{\ga_N (H_{0,\Om}^D-zI_{\Om})^{-1}M_v}  \no \\
& \hspace*{4.3cm} \times
\big[I_{\Om}+\ol{M_u(H_{0,\Om}^D-zI_{\Om})^{-1}M_v}\big]^{-1}\Big)
\no \\
&\quad\quad \times \exp\Big(\tr\Big(\big[\ol{\ga_D
(H_{0,\Om}^N-\ol{z}I_{\Om})^{-1}M_{\ol{u}}}\big]^*
\ol{\ga_N(H_{0,\Om}^D-zI_{\Om})^{-1}M_v}  \lb{6.45}
\\
&\quad\quad \times
\ol{M_u(H_{0,\Om}^D-zI_{\Om})^{-1}M_v}\big[I_{\Om}+\ol{M_u(H_{0,\Om}^D
-zI_{\Om})^{-1}M_v}\big]^{-1}\Big)\Big).  \no
\end{align}
Then, utilizing Corollary \ref{c6.4} with $p_1=12$ and $p_2=12/5$,
one finds,
\begin{align}
\ol{\ga_N(H_{0,\Om}^D-zI_{\Om})^{-1}M_v}
&\in\cB_{12}\big(L^2(\Om;d^nx),L^2(\dOm;d^{n-1}\si)\big),
\\
\ol{\ga_D (H_{0,\Om}^N-\ol{z}I_{\Om})^{-1}M_{\ol{u}}}
&\in\cB_{12/5}\big(L^2(\Om;d^nx),L^2(\dOm;d^{n-1}\si)\big),
\end{align}
and hence using the fact that,
\begin{align}
\big[I_{\Om}+\ol{M_u(H_{0,\Om}^D-zI_{\Om})^{-1}M_v}\big]^{-1} \in
\cB\big(L^2(\Om;d^nx)\big), \quad z\in\bbC\backslash
\big(\si(H_{\Om}^D)\cup\si(H_{0,\Om}^D)\big),
\end{align}
one rearranges the terms in \eqref{6.45} as follows,
\begin{align}
&\frac{\det{}_2\big(I_{\Om}+\ol{M_u(H_{0,\Om}^N-zI_{\Om})^{-1}M_v}\big)}
{\det{}_2\big(I_{\Om}+\ol{M_u(H_{0,\Om}^D-zI_{\Om})^{-1}M_v}\big)}  \no
\\
&\quad = \det{}_2\Big(I_{\dOm} - \ol{\ga_N (H_{0,\Om}^D-zI_{\Om})^{-1}M_v}
\big[I_{\Om}+\ol{M_u(H_{0,\Om}^D-zI_{\Om})^{-1}M_v}\big]^{-1}  \no
\\
& \hspace*{7cm} \times \big[\ol{\ga_D
(H_{0,\Om}^N-\ol{z}I_{\Om})^{-1}M_{\ol{u}}}\,\big]^*\Big)  \no
\\
&\quad\quad \times \exp\Big(\tr\Big(\ol{\ga_N
(H_{0,\Om}^D-zI_{\Om})^{-1} M_v} \;
\ol{M_u(H_{0,\Om}^D-zI_{\Om})^{-1} M_v}  \no
\\
&\quad\quad \times \big[I_{\Om}+ \ol{M_u(H_{0,\Om}^D-zI_{\Om})^{-1}
M_v}\big]^{-1} \big[\ol{\ga_D
(H_{0,\Om}^N-\ol{z}I_{\Om})^{-1}M_{\ol{u}}}\,\big]^*\Big)\Big)  \no
\\
&\quad = \det{}_2\Big(I_{\dOm} - \ol{\ga_N(H_{0,\Om}^D
-zI_{\Om})^{-1}M_{\wti v}}
\big[I_{\Om}+\ol{M_{\wti u}(H_{0,\Om}^D-zI_{\Om})^{-1}M_{\wti
v}}\big]^{-1}   \no
\\
& \hspace*{7cm} \times \big[\ol{\ga_D
(H_{0,\Om}^N-\ol{z}I_{\Om})^{-1}M_{\ol{\wti u}}}\,\big]^*\Big)  \no
\\
&\quad\quad \times
\exp\Big(\tr\Big(\ol{\ga_N(H_{0,\Om}^D-zI_{\Om})^{-1}M_{\wti v}} \;
\ol{M_{\wti u}(H_{0,\Om}^D-zI_{\Om})^{-1}M_{\wti v}}  \lb{6.49}
\\
&\quad\quad \times \big[I_{\Om}+\ol{M_{\wti
u}(H_{0,\Om}^D-zI_{\Om})^{-1}M_{\wti v}}\big]^{-1} \big[\ol{\ga_D
(H_{0,\Om}^N-\ol{z}I_{\Om})^{-1}M_{\ol{\wti u}}}\,\big]^*\Big)\Big).
\no
\end{align}
In the last equality we employed the following simple identities,
\begin{align}
& M_V = M_u M_v = M_{\wti u}  M_{\wti v},  \\
& M_v\big[I_{\Om}+\ol{M_u(H_{0,\Om}^D-zI_{\Om})^{-1}M_v}\big]^{-1}M_u
= M_{\wti v}
\big[I+\ol{M_{\wti u}(H_{0,\Om}^D-zI_{\Om})^{-1}M_{\wti v}}\big]^{-1}
M_{\wti u}.
\end{align}
Utilizing \eqref{6.49} and the following analog of formula \eqref{2.20},
\begin{align}
\ol{(H_{0,\Om}^D-zI_{\Om})^{-1}M_{\wti v}} \big[I_{\Om}+\ol{M_{\wti
u}(H_{0,\Om}^D-zI_{\Om})^{-1}M_{\wti v}}\big]^{-1} =
\ol{(H_{\Om}^D-zI_{\Om})^{-1}M_{\wti v}},
\end{align}
one arrives at \eqref{6.31}, subject to the extra assumption $V\in
L^2(\Om;d^n x)\cap L^6(\Om;d^n x)$.

Finally, assuming only $V\in L^2(\Om;d^n x)$ and utilizing Lemma
\ref{l6.2} and Corollary \ref{c6.4} once again, one obtains
\begin{align}
M_{\wti v}(H_{0,\Om}^D-zI_{\Om})^{-1/6} &\in
\cB_{12}\big(L^2(\Om;d^nx)\big), \lb{6.91a}
\\
M_{\wti u}(H_{0,\Om}^D-zI_{\Om})^{-5/6} &\in
\cB_{12/5}\big(L^2(\Om;d^nx)\big), \lb{6.91b}
\\
\ol{\ga_N(H_{0,\Om}^D-zI_{\Om})^{-1}M_{\wti v}} &\in
\cB_{12}\big(L^2(\Om;d^nx),L^2(\dOm;d^{n-1}\si)\big), \lb{6.91c}
\\
\ol{\ga_D (H_{0,\Om}^N-zI_{\Om})^{-1}M_{\wti u}} &\in
\cB_{12/5}\big(L^2(\Om;d^nx),L^2(\dOm;d^{n-1}\si)\big), \lb{6.91d}
\end{align}
and hence
\begin{equation}
\ol{M_{\wti u}(H_{0,\Om}^D-zI_{\Om})^{-1}M_{\wti v}} \in
\cB_2(L^2(\Om;d^nx)).  \lb{6.92}
\end{equation}
Relations \eqref{6.91a}--\eqref{6.92} prove \eqref{6.70a} and
\eqref{6.70b}. Moreover, since
\begin{align}
\big[I_{\Om}+\ol{M_{\wti u}(H_{0,\Om}^D-zI_{\Om})^{-1}
M_{\wti v}}\big]^{-1} \in
\cB\big(L^2(\Om;d^nx)\big), \quad z\in\bbC\backslash
\big(\si(H_{\Om}^D)\cup\si(H_{0,\Om}^D)\big),
\end{align}
the left- and the right-hand sides of \eqref{6.49}, and hence of
\eqref{6.31}, are well-defined for $V\in L^2(\Om;d^nx)$. Thus,
using \eqref{6.12}, \eqref{6.25}, \eqref{6.26}, the continuity of
$\det{}_2(\cdot)$ with respect to the Hilbert--Schmidt norm
$\|\cdot\|_{\cB_2(L^2(\Om;d^nx))}$, the continuity of
$\tr(\cdot)$ with respect to the trace norm
$\|\cdot\|_{\cB_1(L^2(\Om;d^nx))}$, and an approximation of $V\in
L^2(\Om;d^nx)$ by a sequence of potentials $V_k \in
L^2(\Om;d^nx)\cap L^6(\Om;d^nx)$, $k\in\bbN$, in the norm of
$L^2(\Om;d^nx)$ as $k\uparrow\infty$, then extends the result from
$V\in L^2(\Om;d^nx)\cap L^6(\Om;d^nx)$ to $V\in L^2(\Om;d^nx)$, $n=2,3$.
\end{proof}
%%%%%%%%%%%%%%%%%%%%%%%%%%%%%%%%%%%%%%%

%%%%%%%%%%%%%%%%%%%%%%%%%%%%%%%%%%%%%%% 
\begin{remark}  \lb{r6.12}
Thus, a comparison of Theorem \ref{t6.5} with the one-dimensional case in
Theorem \ref{t6.9} shows that the reduction of Fredholm determinants
associated with operators in $L^2((0,\infty);dx)$ to simple Wronski
determinants, and hence to Jost functions as first observed by Jost and
Pais \cite{JP51}, can be properly extended to higher dimensions and
results in a reduction of appropriate ratios of Fredholm determinants
associated with operators in $L^2(\Om;d^nx)$ to an appropriate
Fredholm determinant associated with an operator in
$L^2(\partial\Om;d^{n-1}\sigma)$. \hfill $\diamond$
\end{remark}
%%%%%%%%%%%%%%%%%%%%%%%%%%%%%%%%%%%%%%%

%%%%%%%%%%%%%%%%%%%%%%%%%%%%%%%%%%%%%%%
\begin{remark} \lb{r6.13}
As in Theorem \ref{t6.5} we assume Hypothesis \ref{h6.1} and suppose
$z\in\bbC\backslash \big(\si(H_{\Om}^D)\cup
\si(H_{0,\Om}^D) \cup \si(H_{0,\Om}^N)\big)$. First we note that
\begin{align}
&\big[(H_{0,\Om}^D-zI_{\Om})^{-1/2}
(H_{\Om}^D-zI_{\Om})(H_{0,\Om}^D-zI_{\Om})^{-1/2} - I_\Om\big]
\in\cB_2\big(L^2(\Om;d^nx)\big),
\\
&\big[(H_{0,\Om}^N-zI_{\Om})^{-1/2}
(H_{\Om}^N-zI_{\Om})(H_{0,\Om}^N-zI_{\Om})^{-1/2} - I_\Om\big]
\in\cB_2\big(L^2(\Om;d^nx)\big).
\end{align}
Indeed, by \eqref{6.6} and \eqref{6.7}, one obtains
\begin{align}
&(H_{0,\Om}^D-zI_{\Om})^{-1/2}
(H_{\Om}^D-zI_{\Om})(H_{0,\Om}^D-zI_{\Om})^{-1/2} - I_\Om  \no
\\
&\quad =
(H_{0,\Om}^D-zI_{\Om})^{-1/2}M_V(H_{0,\Om}^D-zI_{\Om})^{-1/2}
\in\cB_2\big(L^2(\Om;d^nx)\big),
\\
&(H_{0,\Om}^N-zI_{\Om})^{-1/2}
(H_{\Om}^N-zI_{\Om})(H_{0,\Om}^N-zI_{\Om})^{-1/2} - I_\Om  \no
\\
&\quad =
(H_{0,\Om}^N-zI_{\Om})^{-1/2}M_V(H_{0,\Om}^N-zI_{\Om})^{-1/2}
\in\cB_2\big(L^2(\Om;d^nx)\big).
\end{align}
Thus, using \eqref{6.6}--\eqref{6.8} and the cyclicity of
$\det{}_2(\cdot)$, one rearranges the left-hand side of \eqref{6.31} as
follows,
\begin{align}
& \frac{\det{}_2\big(I_{\Om}+\ol{M_u(H_{0,\Om}^N-zI_{\Om})^{-1}M_v}\big)}
{\det{}_2\big(I_{\Om}+\ol{M_u(H_{0,\Om}^D-zI_{\Om})^{-1}M_v}\big)}  \no \\
& \quad =
\frac{\det{}_2\big(I_{\Om}+(H_{0,\Om}^N-zI_{\Om})^{-1/2}M_V
(H_{0,\Om}^N-zI_{\Om})^{-1/2}\big)}
{\det{}_2\big(I_{\Om}+(H_{0,\Om}^D-zI_{\Om})^{-1/2}M_V
(H_{0,\Om}^D-zI_{\Om})^{-1/2}\big)}  \no \\
& \quad =\frac{\det{}_2\big((H_{0,\Om}^N-zI_{\Om})^{-1/2}(H_{\Om}^N
-zI_{\Om})(H_{0,\Om}^N-zI_{\Om})^{-1/2}\big)}
{\det{}_2\big((H_{0,\Om}^D-zI_{\Om})^{-1/2}(H_{\Om}^D-zI_{\Om})
(H_{0,\Om}^D-zI_{\Om})^{-1/2}\big)}.  \lb{6.38}
\end{align}
Again \eqref{6.38} illustrates that symmetrized perturbation determinants
underly Theorem \ref{t6.5}. \hfill $\diamond$
\end{remark}
%%%%%%%%%%%%%%%%%%%%%%%%%%%%%%%%%%%%%%%

%%%%%%%%%%%%%%%%%%%%%%%%%%%%%%%%%%%%%%%
\begin{remark}  \lb{r6.14}
The following observation yields a simple
application of formula \eqref{6.31}. Since by Theorem \ref{t2.6}, for any
$z\in\bbC\backslash \big(\si(H_{\Om}^D)\cup
\si(H_{0,\Om}^D) \cup \si(H_{0,\Om}^N)\big)$,
one has $z\in\si(H_{\Om}^N)$ if and
only if $\det{}_2\big(I_{\Om}+
\ol{M_u(H_{0,\Om}^N -zI_{\Om})^{-1}M_v}\big)=0$, it follows
from \eqref{6.31} that
\begin{align}
\begin{split}
&\text{for all } \, z\in\bbC\backslash \big(\si(H_{\Om}^D)\cup
\si(H_{0,\Om}^D) \cup \si(H_{0,\Om}^N)\big), \text{ one has }
z\in\si(H_{\Om}^N)  \\
&\quad \text{if and only if } \, \det{}_2\big( I_{\dOm} -
\ol{\ga_N(H_{\Om}^D-zI_{\Om})^{-1}M_V
\big[\ga_D(H_{0,\Om}^N-\ol{z}I_{\Om})^{-1}\big]^*}
\big)=0.
\end{split}
\end{align}
One can also prove the following analog of \eqref{6.31}:
\begin{align}
&\frac{\det{}_2\big(I_{\Om}+\ol{M_u(H_{0,\Om}^D-zI_{\Om})^{-1}M_v}\big)}
{\det{}_2\big(I_{\Om}+\ol{M_u(H_{0,\Om}^N-zI_{\Om})^{-1}M_v}\big)} \no \\
&\quad = \det{}_2\big(I_{\dOm} +
\ol{\ga_N(H_{0,\Om}^D-zI_{\Om})^{-1}M_V
\big[\ga_D((H_{\Om}^N-z I_{\Om})^{-1})^*\big]^*}\big)
\lb{6.90} \\
&\quad\quad \times \exp\big(-\tr\big(
\ol{\ga_N(H_{0,\Om}^D-zI_{\Om})^{-1}M_V
(H_{\Om}^N-zI_{\Om})^{-1}M_V
\big[\ga_D(H_{0,\Om}^N-\ol{z}I_{\Om})^{-1}\big]^*}
\big)\big). \no
\end{align}
Then, proceeding as before, one obtains
\begin{align}
& \text{for all } \, z \in \bbC \backslash
\big(\si(H_{\Om}^N)\cup\si(H_{0,\Om}^N)\cup\si(H_{0,\Om}^D)\big),
\text{ one has } z\in\si(H_{\Om}^D) \\
&\quad\text{if and only if}\quad \det{}_2\big(I_{\dOm} +
\ol{\ga_N(H_{0,\Om}^D-zI_{\Om})^{-1}M_V
\big[\ga_D((H_{\Om}^N-z I_{\Om})^{-1})^*\big]^*} \big)=0.  \no
\end{align}
${}$ \hfill $\diamond$
\end{remark}
%%%%%%%%%%%%%%%%%%%%%%%%%%%%%%%%%%%%%

%%%%%%%%%%%%%%%%%%%%%%%%%%%%%%%%%%%%%%%
%%%%%%%%%%%%%%%%%%%%%%%%%%%%%%%%%%%%%%%
\section{An Application to Scattering Theory} \lb{s7}
%%%%%%%%%%%%%%%%%%%%%%%%%%%%%%%%%%%%%%%
%%%%%%%%%%%%%%%%%%%%%%%%%%%%%%%%%%%%%%%

In this section we relate Krein's spectral shift function and hence the
determinant of the scattering operator in connection with quantum
mechanical scattering theory in dimensions $n=2,3$ with appropriate
modified Fredholm determinants.

The results of this section are not new, they were first derived for
$n=3$ by Newton \cite{Ne77} and subsequently for $n=2$ by Cheney
\cite{Ch84}. However, since our method of proof nicely illustrates
the use of infinite determinants in connection with scattering theory and
is different from that in \cite{Ne77} and \cite{Ch84}, and
moreover, since our derivation in the case $n=3$ is performed under
slightly more general hypotheses than in \cite{Ne77}, we thought it
worthwhile to include it at this point.

%%%%%%%%%%%%%%%%%%%%%%%%%%%%%%%%%%%%%%%
\begin{hypothesis} \lb{h7.1}
Fix $\delta>0$. Suppose $V\in \cR_{2,\delta}$ for $n=2$ and $V\in
L^1(\bbR^3;d^3 x)\cap\cR_3$ for $n=3$, where
\begin{align}
\cR_{2,\delta}&=\big\{V\colon \bbR^2\to\bbR \;\text{measurable}\,\big|\,
V^{1+\delta}, (1+|\cdot|^{\delta})V \in L^1(\bbR^2; d^2x)\big\}, \\
\cR_3&=\bigg\{V\colon \bbR^3\to\bbR\; \text{measurable}\,\bigg|\,
\int_{\bbR^6}  d^3 x d^3 x' \,
|V(x)||V(x')||x-x'|^{-2}<\infty\bigg\}.
\end{align}
\end{hypothesis}
%%%%%%%%%%%%%%%%%%%%%%%%%%%%%%%%%%%%%%

We introduce $H_0$ as the following nonnegative self-adjoint operator in
the Hilbert space $L^2(\bbR^n;d^nx)$,
\begin{equation}
H_0 = -\Delta, \quad \dom(H_0) = H^{2}(\bbR^n), \quad n=2,3.
\end{equation}
Moreover, let $A=M_u$ and $B=B^*=M_v$ denote the operators of
multiplication by
$u=\sgn(V)\abs{V}^{1/2}$ and $v=\abs{V}^{1/2}$ in
$L^2(\bbR^n;d^nx)$, respectively, so that $M_V=B A=M_u M_v$.  Then, (cf.\
\cite[Theorem I.21]{Si71} for
$n=3$ and
\cite{Si76} for $n=2$),
\begin{align}
\dom(A)=\dom(B) \supseteq H^{1}(\bbR^n)\supset \dom(H_0), \lb{7.3A}
\end{align}
and hence, Hypothesis \ref{2.1}\,$(i)$ is satisfied for $H_0$. It
follows from Hypothesis \ref{h7.1} that
\begin{align}
&\ol{M_u(H_0-zI)^{-1}M_v}\in\cB_2\big(L^2(\bbR^n;d^nx)\big),  \quad
z\in\bbC\backslash [0,\infty),  \lb{7.4}
\end{align}
where $I$ now denotes the identity operator in $L^2(\bbR^n;d^nx)$, and
hence, Hypothesis \ref{2.1}\,$(ii)$ is satisfied. Taking
$z\in\bbC\backslash[0,\infty)$ with a sufficiently large absolute
value, one also verifies Hypothesis \ref{2.1}\,$(iii)$. Thus, applying
Theorem \ref{2.3} and Remark \ref{r2.4}\,$(i)$, one obtains a
self-adjoint  operator $H$ (which is an extension of $H_0+V$ on
$\dom(H_0)\cap\dom(V)$).

%%%%%%%%%%%%%%%%%%%%%%%%%%%%%%%%%%%%%
\begin{theorem} \lb{t7.2}
Assume Hypothesis \ref{h7.1} and let $z\in\bbC\backslash\si(H)$ and
$n=2,3$. Then,
\begin{align}
(H-zI)^{-1} - (H_0-zI)^{-1} \in\cB_1\big(L^2(\bbR^n;d^nx)\big),  \lb{7.5}
\end{align}
and there is a unique real-valued spectral shift function
\begin{equation}
\xi(\cdot,H,H_0)\in L^1\big(\bbR; (1+\la^2)^{-1}d\la)   \lb{7.6}
\end{equation}
such that $\xi(\la,H,H_0)=0$ for $\la< \inf (\sigma(H))$, and
\begin{align}
\tr\big((H-zI)^{-1} - (H_0-zI)^{-1}\big) = - \int_{\sigma(H)}
\frac{d\la\,\xi(\la,H,H_0)}{(\la-z)^2}.   \lb{7.7}
\end{align}
\end{theorem}
%%%%%%%%%%%%%%%%%%%%%%%%%%%%%%%%%%%%%%

We recall that $\xi(\cdot,H,H_0)$ is called the spectral shift function
for the pair of self-adjoint operators $(H,H_0)$. For background
information on $\xi(\cdot,H,H_0)$ and its connection with the scattering
operator at fixed energy, we refer, for instance, to \cite[Sect.\
19.1]{BW83}, \cite{BK62}, \cite{BY93}, \cite[Ch.\ 8]{Ya92}.

%%%%%%%%%%%%%%%%%%%%%%%%%%%%%%%%%%%%%%
\begin{lemma}  \lb{l7.3}
Assume Hypothesis \ref{h7.1} and let $z\in\bbC\backslash\si(H)$ and
$n=2,3$. Then,
\begin{align}
 \ol{M_u(H_0-zI)^{-1}M_v} &\in\cB_2\big(L^2(\bbR^n;d^nx)\big),  \lb{7.8}
\\
 (H_0-zI)^{-1}M_V(H_0-zI)^{-1}
& \in\cB_1\big(L^2(\bbR^n;d^nx)\big),   \lb{7.9}
\end{align}
and
\begin{align}
\begin{split}
&\frac{d}{dz}\ln\big(\det{}_2\big(I+\ol{M_u(H_0-zI)^{-1}M_v}\big)\big)
\\
&\quad = -\tr\big((H-zI)^{-1} - (H_0-zI)^{-1} +
(H_0-zI)^{-1}M_V(H_0-zI)^{-1}\big).   \lb{7.10}
\end{split}
\end{align}
\end{lemma}
%%%%%%%%%%%%%%%%%%%%%%%%%%%%%%%%%%%%%%

The key ingredient in proving \eqref{7.5} is the fact that
\begin{equation}
M_u(H_0-zI)^{-1}, \, \ol{(H_0-zI)^{-1}M_v} \in\cB_2\big(L^2(\bbR^n;d^n
x)\big),
\quad z\in\bbC\backslash [0,\infty), \; n=2,3.  \lb{7.10a}
\end{equation}
This follows from either \cite[Theorem 4.1]{Si79} (or \cite[Theorem
XI.20]{RS79}), or explicitly by an inspection of the corresponding integral
kernels. For instance, the one for $M_u(H_0-zI)^{-1}$ reads:
\begin{align}
& \big(M_u(H_0-zI)^{-1}\big)(x,x')= \begin{cases}
u(x)(i/4)H_0^{(1)}(z^{1/2}|x-x'|), & x \neq x', \; x,
x'\in\bbR^2,  \\[2mm]
u(x)e^{i z^{1/2}|x-x'|}/[4\pi |x-x'|], &x \neq x', \; x, x'\in\bbR^3,
\end{cases} \no \\
& \hspace*{7cm}  z\in\bbC\backslash [0,\infty), \; \Im(z^{1/2})>0,
\lb{7.10b}
\end{align}
where $H_0^{(1)}(\cdot)$ denotes the Hankel function of order zero and
first kind (see, e.g., \cite[Sect.\ 9.1]{AS72}). Hence, one only needs to
apply equation \eqref{2.13} to conclude \eqref{7.5} and hence \eqref{7.9}
(by factoring $M_V=M_u M_v$). (We note that \eqref{7.5} is proved in
\cite[Sect.\ XI.6]{RS79} and \cite[Theorem II.37]{Si71} for $n=3$.)
Relation \eqref{7.8} is then clear from
$V\in\cR_3$ for $n=3$ and follows from \cite{Si76} for $n=2$. Equation
\eqref{7.10} is discussed in \cite{BGGSS87} for $n=2,3$. The trace
formula \eqref{7.7} is a celebrated result of Krein
\cite{Kr53}, \cite{Kr62}; detailed accounts of it can be found in
\cite[Sect.\ 19.1.5]{BW83}, \cite{BY93}, \cite{Kr83}, \cite[Ch.\ 8]{Ya92}.

%%%%%%%%%%%%%%%%%%%%%%%%%%%%%%%%%%%%%%
\begin{lemma} \lb{l7.4}
Assume Hypothesis \ref{h7.1}. Then the following formula holds for a.e.\
$\lambda\in\bbR$,
\begin{align}
\begin{split}
2\pi i\xi(\la,H,H_0) =& \; \ln\left(
\frac{\det{}_2\big(I+\ol{M_u(H_0-(\la+i0)I)^{-1}M_v})}
{\det{}_2\big(I+\ol{M_u(H_0-(\la-i0)I)^{-1}M_v})} \right) \lb{7.11}
\\
&+ \frac{i}{2\pi}\int_{\bbR^n} d^nx\,V(x) \times \begin{cases} \pi,
& \la>0, \; n=2,\\
\la^{1/2}, & \la>0, \; n=3,\\
0, & \la\leq0, \; n=2,3.
\end{cases}
\end{split}
\end{align}
\end{lemma}
%%%%%%%%%%%%%%%%%%%%%%%%%%%%%%%%%%%%%
\begin{proof}
It follows from Theorem \ref{t7.2} and Lemma \ref{l7.3}, that for
$z\in\bbC\backslash\si(H)$,
\begin{align}
\begin{split}
\int_\bbR \frac{d\la\,\xi(\la,H,H_0)}{(\la-z)^2} =&
\frac{d}{dz}\ln\big(\det{}_2\big(I+\ol{M_u(H_0-zI)^{-1}M_v}\big)\big)
\lb{7.12}
\\
&+ \tr\big((H_0-zI)^{-1}M_V(H_0-zI)^{-1}\big).
\end{split}
\end{align}
First, we rewrite the left-hand side of \eqref{7.12}. Since
$\xi(\cdot,H,H_0)\in L^1\big(\bbR;\frac{d\la}{1+\la^2}\big)$, one
has the following formula,
\begin{align}
\int_\bbR \frac{d\la\,\xi(\la,H,H_0)}{(\la-z)^2} = \frac{d}{dz}
\int_\bbR
d\la\,\xi(\la,H,H_0)\left(\frac{1}{\la-z}-\frac{\la}{1+\la^2}\right),
\quad z\in\bbC\backslash\si(H). \lb{7.13}
\end{align}

Next, we compute the second term on the right-hand side of
\eqref{7.12}. By \eqref{7.10a} and the cyclicity of the trace,
\begin{align}
\tr\big((H_0-zI)^{-1}M_V(H_0-zI)^{-1}\big) =
\tr\big(\ol{M_u(H_0-zI)^{-2}M_v}\big), \quad z\in\bbC\backslash[0,\infty).
\end{align}
Then $\ol{M_u(H_0-zI)^{-2}M_v}$ = $\ol{M_u\frac{d}{dz}(H_0-zI)^{-1}M_v}$
has the integral kernel
\begin{align}
& \big(\ol{M_u(H_0-zI)^{-2}M_v}\big)(x,x')= \begin{cases}
u(x)\frac{i{H_0^{(1)}}'(z^{1/2}\abs{x-x'})\abs{x-x'}}{8 z^{1/2}}v(x'),
&x, x'\in\bbR^2, \\[2mm]
 u(x)\frac{i\exp(i z^{1/2}\abs{x-x'})}{8\pi z^{1/2}}v(x'), &
 x, x'\in\bbR^3, \end{cases}  \no \\
& \hspace*{5.8cm} x\neq x', \;
 z\in\bbC\backslash [0,\infty), \; \Im(z^{1/2})>0,
\end{align}
and hence, utilizing \cite[p.\ 1086]{DS88}, one computes for
$z\in\bbC\backslash[0,\infty)$,
\begin{align}
\begin{split}
\tr\big((H_0-zI)^{-1}M_V(H_0-zI)^{-1}\big) &=
\frac{1}{4\pi}\int_{\bbR^n} d^nx\,V(x)
\times \begin{cases} -z^{-1}, & n=2
\\
i(2 z^{1/2})^{-1}, & n=3
\end{cases}
\\
&= \frac{1}{4\pi}\int_{\bbR^n} d^nx\,V(x) \times \frac{d}{dz}
\begin{cases} -\ln(z), & n=2, \\ i z^{1/2}, & n=3. \end{cases}
\lb{7.17}
\end{split}
\end{align}

Finally, using \eqref{7.12}, \eqref{7.13}, and \eqref{7.17},
one obtains for $z\in\bbC\backslash\si(H)$,
\begin{align}
\begin{split}
& \int_\bbR d\la\, \xi(\la,H,H_0)
\left(\frac{1}{\la-z}-\frac{\la}{1+\la^2}\right) + C    \\
& \quad =\ln\big(\det{}_2\big(I+\ol{M_u(H_0-zI)^{-1}M_v}\big)\big)
+\frac{1}{4\pi}\int_{\bbR^n} d^nx\,V(x)
\times \begin{cases} -\ln(z), & n=2, \\
i z^{1/2}, & n=3,
\end{cases}  \lb{7.21}
\end{split}
\end{align}
where $C\in\bbC$ denotes an appropriate constant. To complete the proof
we digress for a moment and recall the Stieltjes inversion formula for
Herglotz functions $m$ (i.e., analytic maps $m\colon \bbC_+\to\bbC_+$,
where $\bbC_+$ denotes the open complex upper half-plane). Such functions
$m$ permit the Nevanlinna, respectively, Riesz-Herglotz
representation
\begin{align}
\begin{split}
&m(z)=c+dz+\int_{{\mathbb{R}}}
d\omega (\lambda) \bigg(\frac{1}{\lambda -z}-\frac{\lambda}
{1+\lambda^2}\bigg), \quad z\in\bbC_+, \lb{7.22} \\[2mm]
& \, c=\Re[m(i)],\quad d=\lim_{\eta \uparrow
\infty}m(i\eta )/(i\eta ) \geq 0,
\end{split}
\end{align}
with a nonnegative measure $d\omega$ on ${\mathbb{R}}$ satisfying
\begin{equation} \lb{7.23}
\int_{{\mathbb{R}}} \frac{d\omega (\lambda )}{1+\lambda^2}<\infty.
\end{equation}
The absolutely continuous part $d\omega_{ac}$ of
$d\omega$ with respect to Lebesgue measure $d\lambda$ on
${\mathbb{R}}$ is then known to be given by
\begin{equation}\lb{7.24}
d\omega_{ac}(\lambda)=\pi^{-1}\Im[m(\lambda+i0)]\,d\lambda.
\end{equation}
In addition, one extends $m$ to the open lower complex half-plane $\bbC_-$
by
\begin{equation}
m(z)=\ol{m(\ol z)}, \quad z\in\bbC_-.   \lb{7.25}
\end{equation}
(We refer, e.g., to \cite[Sect.\ 69]{AG81} for details on
\eqref{7.22}--\eqref{7.25}.) Thus, in order to apply
\eqref{7.22}--\eqref{7.25} to the computation of
$\xi(\cdot,H,H_0)$ in \eqref{7.21} it suffices to decompose
$\xi(\cdot,H,H_0)=\xi_+(\cdot,H,H_0)-\xi_-(\cdot,H,H_0)$ into its
positive and negative parts $\xi_{\pm}(\cdot,H,H_0)\geq 0$ and
separately consider the absolutely continuous measures
$\xi_{\pm}(\cdot,H,H_0) d\lambda$.  Thus, letting
$z=\la\pm i\eps$, taking the limit $\eps\downarrow 0$ in \eqref{7.21},
and subtracting the corresponding results, yields \eqref{7.11}.
\end{proof}
%%%%%%%%%%%%%%%%%%%%%%%%%%%%%%%%%%%%%%

We conclude with the following result:

%%%%%%%%%%%%%%%%%%%%%%%%%%%%%%%%%%%%%%
\begin{corollary} \lb{c7.6}
Assume Hypothesis \ref{h7.1}. Then, for a.e. $\la>0$,
\begin{align}
\begin{split}
\det(S(\la)) =& \; \frac{\det{}_2\big(I+\ol{M_u(H_0-(\la-i0)I)^{-1}
M_v}\big)}
{\det{}_2\big(I+\ol{M_u(H_0-(\la+i0)I)^{-1}M_v}\big)}  \lb{7.19}
\\
&\times \begin{cases} \exp\big(-\frac{i}{2}\int_{\bbR^n}
d^nx\,V(x)\big), & n=2,
\\[2mm]
\exp\big(-\frac{i \la^{1/2}}{2\pi}\int_{\bbR^n} d^nx\,V(x)\big), &
n=3.
\end{cases}
\end{split}
\end{align}
\end{corollary}
%%%%%%%%%%%%%%%%%%%%%%%%%%%%%%%%%%%%%%
\begin{proof}
Hypothesis \ref{h7.1} implies that the scattering operator $S(\la)$ at
fixed energy $\lambda>0$ in $L^2(S^{n-1};d^{n-1}\omega)$ satisfies
\begin{align}
[S(\la)-I] &\in\cB_1\big(L^2(S^{n-1};d^{n-1}\omega)\big) \, \text{ for
a.e.\ $\lambda>0$}
\end{align}
and
\begin{align}
\det(S(\la)) = \exp(-2\pi i \xi(\la,H,H_0)) \, \text{ for a.e.\ $\la>0$}
\lb{7.27}
\end{align}
(cf., e.g., \cite[Sects.\ 19.1.4, 19.1.5]{BW83}, \cite{BK62}, \cite{BY93},
\cite[Ch.\ 8]{Ya92}), where $S^{n-1}$ denotes the unit sphere in $\bbR^n$
and $d^{n-1}\omega$ the corresponding surface measure on $S^{n-1}$.
Relation \eqref{7.19} then follows from Lemma \ref{l7.4} and \eqref{7.27}.
\end{proof}
%%%%%%%%%%%%%%%%%%%%%%%%%%%%%%%%%%%%%%
We note again that Corollary 7.5 was derived earlier using
different means by Cheney \cite{Ch84} for $n=2$ and by Newton \cite{Ne77}
for $n=3$. (The stronger conditions $V\in L^2(\bbR^3; dx^3)$ and
the existence of $a>0$ and $0<C<\infty$ such that for all $y\in\bbR^3$,
$\int_{\bbR^3} d^3x \, |V(x)| [(|x|+|y|+a)/(|x-y|)]^2 \leq C$, are assumed
in \cite{Ne77}.)

%%%%%%%%%%%%%%%%%%%%%%%%%%%%%%%%%%%%%%%%
%%%%%%%%%%%%%%%% appendices %%%%%%%%%%%%%%%%% 
\appendix
%%%%%%%%%%%%%%%% appendix A %%%%%%%%%%%%%%%%%
\section{Properties of the Dirichlet and Neumann Laplacians}
\lb{sA}
\renewcommand{\theequation}{A.\arabic{equation}}
\renewcommand{\thetheorem}{A.\arabic{theorem}}
\setcounter{theorem}{0} \setcounter{equation}{0}
%%%%%%%%%%%%%%%%%%%%%%%%%%%%%%%%%%%%%%%%
%%%%%%%%%%%%%%%%%%%%%%%%%%%%%%%%%%%%%%%%

The purpose of this appendix is to derive some basic domain properties of
Dirichlet and Neumann Laplacians on $C^{1,r}$-domains $\Om\subset\bbR^n$
and to prove Lemma \ref{l6.1}. Throughout this appendix we assume
$n\geq 2$, but we note that $n$ is restricted to $n=2,3$ in
Sections \ref{s6} and \ref{s7}.

In this manuscript we use the following notation for the standard
Sobolev Hilbert spaces ($s\in\bbR$),
\begin{align}
H^{s}(\bbR^n) &=\bigg\{U\in \cS^\prime (\bbR^n) \,\bigg|\,
\norm{U}_{H^{s}(\bbR^n)}^2 = \int_{\bbR^n} d^n \xi \, \big|\hatt
U(\xi)\big|^2\big(1+\abs{\xi}^{2s}\big) <\infty \bigg\},
\\
H^{s}(\Om) &=\left\{u\in \cD^\prime(\Om) \,|\, u=U|_\Om \text{
for some } U\in H^{s}(\bbR^n) \right\},
\\
H_0^{s}(\Om) &=\{u\in H^s(\bbR^n)\,|\, \supp\,(u)\subseteq\ol{\Om}\}.
\end{align}
Here $\cD^\prime(\Om)$ denotes the usual set of distributions on
$\Omega\subseteq \bbR^n$, $\Omega$ open and nonempty,
$\cS^\prime (\bbR^n)$ is the space of tempered distributions on
$\bbR^n$, and $\hatt U$ denotes the Fourier transform of $U\in
\cS^\prime (\bbR^n)$. It is then immediate that
\begin{equation}\label{incl-xxx}
H^{s_1}(\Omega)\hookrightarrow H^{s_0}(\Omega) \, \text{ for } \,
-\infty<s_0\leq s_1<+\infty,
\end{equation}
continuously and densely. 

Before we present a proof of Lemma \ref{l6.1}, we recall the
definition of a $C^{1,r}$-domain $\Omega\subseteq\bbR^n$, $\Om$ open and
nonempty, for convenience of the reader: Let ${\mathcal N}$ be a space of
real-valued functions in $\bbR^{n-1}$.  One calls a bounded domain
$\Omega\subset\bbR^n$ of class ${\mathcal N}$ if there exists a finite
open covering $\{{\mathcal O}_j\}_{1\leq j\leq N}$ of the boundary
$\partial\Omega$ of $\Om$ with the property that, for every
$j\in\{1,...,N\}$, ${\mathcal O}_j\cap\Omega$ coincides with the portion of
${\mathcal O}_j$ lying in the over-graph of a function
$\varphi_j\in{\mathcal N}$ (considered in a new system
of coordinates obtained from the original one via a rigid motion). Two
special cases are going to play a particularly important role in the
sequel. First, if ${\mathcal  N}$ is ${\rm Lip}\,(\bbR^{n-1})$,
the space of real-valued functions satisfying a (global) Lipschitz
condition in $\bbR^{n-1}$, we shall refer to $\Omega$ as being a Lipschitz
domain; cf.\ \cite[p.\ 189]{St70}, where such domains are called
``minimally smooth''. Second, corresponding to  the case when
${\mathcal N}$ is the subspace of ${\rm Lip}\,(\bbR^{n-1})$ consisting of
functions whose first-order derivatives satisfy a (global) H\"older
condition of order $r\in(0,1)$, we shall say that $\Omega$ is of class
$C^{1,r}$. The classical theorem of  Rademacher of almost everywhere
differentiability of Lipschitz functions ensures that, for any  Lipschitz
domain $\Omega$, the surface measure
$d^{n-1}\sigma$ is well-defined on  $\partial\Omega$ and that there exists an
outward  pointing normal vector $\nu$ at almost every point of
$\partial\Omega$. For a Lipschitz domain $\Omega\subset\bbR^n$ it is
known that
\begin{align}\lb{dual-xxx}
\bigl(H^{s}(\Omega)\bigr)^*=H^{-s}(\Omega), \quad
{\textstyle{-\frac12}}<s<{\textstyle{\frac12}}.
\end{align}
See \cite{Tr02} for this and other related properties.

Next, assume that $\Omega\subset\bbR^n$ is the domain lying above the
graph of a function $\varphi:\bbR^{n-1}\to\bbR$ of class $C^{1,r}$.
Then for $0\leq s<1+r$, the Sobolev space $H^s(\partial\Omega)$ consists
of functions $f\in L^2(\partial\Omega;d^{n-1}\sigma)$ such that
$f(x',\varphi(x'))$, as a function of $x'\in\bbR^{n-1}$, belongs
to $H^s(\bbR^{n-1})$. This definition is easily adapted to the
case when $\Omega$ is a domain of class $C^{1,r}$ whose boundary is compact,
by using a smooth partition of unity. Finally, for $-1-r<s<0$, we set
$H^s(\partial\Omega)=(H^{-s}(\partial\Omega))^*$. For additional
background information in this context we refer, for instance, to
\cite[Ch.\ 3]{Mc00}, \cite[Sect.\ I.4.2]{Wl87}.

Assuming Hypothesis \ref{h6.1}\,$(i)$ (i.e., $\Om$ is an open
nonempty $C^{1,r}$-domain for some $(1/2)<r<1$ with compact boundary
$\partial\Om)$, we introduce the Dirichlet and Neumann Laplacians
$\wti H_{0,\Om}^D$ and
$\wti H_{0,\Om}^N$ associated with the domain $\Om$ as the unique
self-adjoint operators on $L^2(\Om;d^nx)$ whose quadratic form equals
$q(f,g)=\int_\Om d^nx\,\ol{\nabla f}\cdot \nabla g$ with (form) domains
$H_0^{1}(\Om)$ and $H^{1}(\Om)$, respectively. Then,
\begin{align}
\dom(\wti H_{0,\Om}^D) &= \{u\in H_0^{1}(\Om) \,|\, \text{there exists} \,
f\in L^2(\Om;d^nx) \text{ such that } \no \\
&\hspace*{1.85cm} q(u,v)=(f,v)_{L^2(\Om;d^nx)} \text{ for all } v\in
H_0^{1}(\Om)\},  \\
\dom(\wti H_{0,\Om}^N) &= \{u\in H^{1}(\Om) \,|\, \text{there exists}\,
f\in L^2(\Om;d^nx) \text{ such that } \no \\
&\hspace*{1.85cm} q(u,v)=(f,v)_{L^2(\Om;d^nx)} \text{ for all } v\in
H^{1}(\Om)\},
\end{align}
with $(\cdot,\cdot)_{L^2(\Om;d^nx)}$ denoting the scalar product in
$L^2(\Om;d^nx)$. Equivalently, we introduce the densely defined closed
linear operators
\begin{equation}
D=\nabla, \; \dom(D)=H_0^{1}(\Om) \, \text{ and } \,
N=\nabla, \;  \dom(N)=H^{1}(\Om)
\end{equation}
from $L^2(\Om;d^nx)$ to $L^2(\Om;d^nx)^n$ and note that
\begin{equation}
\wti H_{0,\Om}^D = D^*D \, \text{ and } \, \wti H_{0,\Om}^N = N^*N.
\end{equation}
For details we refer to \cite[Sects.\ XIII.14, XIII.15]{RS78}.
Moreover, with ${\rm div(\cdot)}$ denoting the divergence operator,
\begin{align}
\dom(D^*)=\{w\in L^2(\Omega;d^nx)^n\,|\,{\rm div} (w)\in
L^2(\Omega;d^nx)\},
\end{align}
and hence,
\begin{align}
\dom(\wti H_{0,\Om}^D) &= \{u\in\dom(D) \,|\,Du\in\dom(D^*)\} \no
\\
&= \{u\in H_0^{1}(\Om) \,|\, \Delta u\in L^{2}(\Om;d^nx)\}.
\lb{domHD}
\end{align}
One can also define the following map
\begin{equation}
\begin{cases} \{w\in
L^2(\Omega;d^nx)^n\,|\,{\rm div}(w)\in (H^1(\Omega))^*\} \to
H^{-1/2}(\partial\Omega)
=\big(H^{1/2}(\partial\Omega)\big)^* \\
\hspace*{5.75cm} w\mapsto \nu\cdot w  \end{cases}  \lb{A.11}
\end{equation}
by setting
\begin{align}
\langle\nu\cdot
w,\phi\rangle=\int_{\Omega}d^nx\,w(x)\cdot\nabla\Phi(x) +
\langle {\rm div}(w)\,,\,\Phi\rangle
\lb{A.11a}
\end{align}
whenever $\phi\in H^{1/2}(\partial\Omega)$ and $\Phi\in
H^{1}(\Omega)$ is such that $\ga_D\Phi=\phi$. The last paring
in (\ref{A.11a}) is in the duality sense (which, in turn, is compatible
with the (bilinear) distributional pairing). It should be remarked
that the above definition is independent of the particular extension
$\Phi\in H^{1}(\Omega)$ of $\phi$. Indeed, by linearity this comes
down to proving that
\begin{equation}\label{ibp}
\langle {\rm div}(w)\,,\,\Phi\rangle
=-\int_{\Omega}d^nx\,w(x)\cdot\nabla\Phi(x)
\end{equation}
if $w\in L^2(\Omega;d^nx)^n$ has ${\rm div}(w)\in \big(H^1(\Omega)\big)^*$
and $\Phi\in H^{1}(\Omega)$ has $\ga_D\Phi=0$. To see this we rely
on the existence of a sequence $\Phi_j\in C^\infty_0(\Omega)$ such
that $\Phi_j\underset{j\uparrow\infty}{\rightarrow} \Phi$ in
$H^{1}(\Omega)$. When $\Omega$ is a bounded Lipschitz domain, this is
well-known (see, e.g., \cite[Remark 2.7]{JK95} for a rather general result
of this nature), and this result is easily extended to the case when
$\Omega$ is an unbounded Lipschitz domain with a compact boundary. For if
$\xi\in C^\infty_0(B(0;2))$ is such that $\xi= 1$ on $B(0;1)$ and
$\xi_j(x)=\xi(x/j)$, $j\in\bbN$ (here $B(x_0;r_0)$ denotes the ball in
$\bbR^n$ centered at $x_0\in\bbR^n$ of radius $r_0>0$), then
$\xi_j\Phi\underset{j\uparrow\infty}{\rightarrow}\Phi$ in
$H^{1}(\Omega)$ and matters are reduced to approximating $\xi_j\Phi$ in
$H^{1}(B(0;2j)\cap\Omega)$ with test functions supported in
$B(0;2j)\cap\Omega$, for each fixed $j\in\bbN$. Since
$\ga_D(\xi_j\Phi)=0$, the result for bounded Lipschitz domains applies.

Returning to the task of proving (\ref{ibp}), it suffices to prove a
similar identity with $\Phi_j$ in place of $\Phi$. This, in turn, follows
from the definition of ${\rm div}(\cdot)$ in the sense of distributions
and the fact that the duality between $(H^1(\Omega))^*$ and $H^1(\Omega)$
is compatible with the duality between distributions and test functions.

Going further, we can introduce a (weak) Neumann trace operator
$\wti\gamma_N$ as follows:
\begin{align}\lb{A.16}
\wti\gamma_N:\{u\in H^{1}(\Om) \,|\, \Delta u \in (H^1(\Om))^*\}\to
H^{-1/2}(\dOm),\quad \wti\gamma_N u=\nu\cdot \nabla u,
\end{align}
with the dot product  understood in the sense of \eqref{A.11}.
We emphasize that the weak Neumann trace operator $\wti\gamma_N$ in
\eqref{A.16} is an extension of the operator $\gamma_N$
introduced in \eqref{6.2}. Indeed, to see that
$\dom(\gamma_N)\subset\dom(\wti\gamma_N)$, we note that
if $u\in H^{s+1}(\Omega)$ for some $1/2<s<3/2$, then
$\Delta u\in H^{-1+s}(\Omega)=\bigl(H^{1-s}(\Omega)\bigr)^*\hookrightarrow
\bigl(H^{1}(\Omega)\bigr)^*$, by (\ref{dual-xxx}) and (\ref{incl-xxx}).
With this in hand, it is then easy to show that
$\wti\gamma_N$ in \eqref{domHN} and $\gamma_N$ in \eqref{6.2} agree
(on the smaller domain), as claimed.

We now return to the mainstream discussion.
From the above preamble it follows that
\begin{equation}
\dom(N^*)=\{w\in L^2(\Omega;d^nx)^n\,|\,{\rm div}(w)\in
L^2(\Omega;d^nx) \mbox{ and }\nu\cdot w=0\}, 
\end{equation}
where the dot product operation is understood in the sense of \eqref{A.11}.
Consequently, with $\wti H_{0,\Om}^N=N^*N$, we have
\begin{align}
\dom(\wti H_{0,\Om}^N) & = \{u\in\dom(N) \,|\,Nu\in\dom(N^*)\} \no
\\
& =\{u\in H^{1}(\Om) \,|\, \Delta u\in L^{2}(\Om;d^nx)\mbox{ and }
\wti\gamma_N u=0\}.
\lb{domHN}
\end{align}

Next, we will prove that $H_{0,\Om}^D = \wti H_{0,\Om}^D$ and
$H_{0,\Om}^N = \wti H_{0,\Om}^N$, where $H_{0,\Om}^D$ and
$H_{0,\Om}^N$ denote the operators introduced in \eqref{6.3} and
\eqref{6.4}, respectively. Since it follows from the first Green's
formula (cf., e.g., \cite[Theorem 4.4]{Mc00}) that $H_{0,\Om}^D \subseteq
\wti H_{0,\Om}^D$ and $H_{0,\Om}^N \subseteq \wti H_{0,\Om}^N$, it
remains to show that $H_{0,\Om}^D \supseteq \wti H_{0,\Om}^D$ and
$H_{0,\Om}^N \supseteq \wti H_{0,\Om}^N$. Moreover, it follows from
comparing \eqref{6.3} with \eqref{domHD} and \eqref{6.4} with
\eqref{domHN}, that one needs only to show that $\dom(\wti
H_{0,\Om}^D)$, $\dom(\wti H_{0,\Om}^N)\subseteq H^{2}(\Omega)$.

%%%%%%%%%%%%%%%%%%%%%%%%%%%%%%%%%%%%%%%
\begin{lemma} \lb{lA.1}
Assume Hypothesis \ref{h6.1}\,$(i)$. Then,
\begin{equation}
\dom(\wti H_{0,\Om}^D)\subseteq H^{2}(\Omega), \quad
 \dom(\wti H_{0,\Om}^N)\subseteq
H^{2}(\Omega).
\end{equation}
In particular,
\begin{equation}
H_{0,\Om}^D = \wti H_{0,\Om}^D, \quad
H_{0,\Om}^N = \wti H_{0,\Om}^N.
\end{equation}
\end{lemma}
%%%%%%%%%%%%%%%%%%%%%%%%%%%%%%%%%%%%%%%
\begin{proof}
Consider $u\in\dom(\wti H_{0,\Om}^N)$ and set $f=\Delta u-u\in
L^2(\Omega;d^nx)$. Viewing $f$ as an element in $\big(H^1(\Omega)\big)^*$,
the classical Lax-Milgram Lemma implies that $u$ is the unique solution of
the boundary-value problem
\begin{equation}\label{BVP}
\left\{
\begin{array}{l}
(\Delta -I_{\Omega})u=f\in L^2(\Omega)\hookrightarrow
\bigl(H^{1}(\Omega)\bigr)^*, \\
u\in H^{1}(\Omega), \\
\wti\gamma_N u=0.
\end{array}
\right.
\end{equation}
One convenient way to show that actually
\begin{equation}\label{goal}
u\in H^{2}(\Omega),
\end{equation}
is to use layer potentials. Specifically, let $E(x)$,
$x\in\bbR^n\backslash\{0\}$, be the fundamental solution of the
Helmholtz operator $\Delta -I_{\Omega}$ in $\bbR^n$ and denote by
$(\Delta-I_{\Omega})^{-1}$ the operator of convolution with $E$. Let us
also define the associated single layer potential
\begin{equation}\label{sing-layer}
{\mathcal
S}g(x)=\int_{\partial\Omega}d^{n-1}\sigma_y\,E(x-y)g(y),\quad
x\in\Omega,
\end{equation}
where $g$ is an arbitrary measurable function on $\partial\Omega$. As
is well-known (the interested reader may consult, e.g., \cite{MMT01},
\cite{Ve84} for jump relations in the context of Lipschitz domains), if
\begin{equation}\label{Ksharp}
K^{\#}g(x)=\int_{\partial\Omega}d^{n-1}\sigma_y\,\partial_{\nu_x}E(x-y)g(y),
\quad x\in\partial\Omega
\end{equation}
stands for the so-called adjoint double layer on $\partial\Omega$,
the following jump formula holds
\begin{equation}\label{jump}
\wti\gamma_N {\mathcal S}g=({\textstyle{\frac12}}I_\dOm+K^{\#})g.
\end{equation}

Now, the solution $u$ of (\ref{BVP}) is given by
\begin{equation}\label{sol}
u=(\Delta-I_{\Omega})^{-1}f-{\mathcal S}g
\end{equation}
for a suitable chosen $g$. In order to continue, we recall that the
classical Calder\'on-Zygmund theory yields that, locally,
$(\Delta-I_{\Omega})^{-1}$ is smoothing of order $2$ on the scale of
Sobolev spaces, and since $E$ has exponential decay at infinity, it
follows that $(\Delta-I_{\Omega})^{-1}f\in H^{2}(\Omega)$ whenever $f\in
L^2(\Omega;d^nx)$. We shall then require that
\begin{equation}
\ga_N{\mathcal S}g=\ga_N (\Delta-I_{\Omega})^{-1}f \mbox{ or }
({\textstyle{\frac12}}I_\dOm+K^{\#})g=h= \ga_N
(\Delta-I_{\Omega})^{-1}f\in H^{1/2}(\partial\Omega).
\end{equation}
Thus, formally, $g=(\frac12 I_\dOm+K^{\#})^{-1}h$ and (\ref{goal})
follows as soon as we prove that
\begin{equation}\label{goal-1}
{\textstyle{\frac12}}I_\dOm+K^{\#}\mbox{ is invertible on }
H^{1/2}(\partial\Omega)
\end{equation}
and that the operator
\begin{equation}\label{goal-2}
{\mathcal S}\colon H^{1/2}(\partial\Omega)\rightarrow
H^{2}(\Omega)
\end{equation}
is well-defined and bounded. That (\ref{goal-1}) holds is essentially
well-known. See, for instance, \cite[Proposition 4.5]{Ta96}
which requires that $\Om$ is of class $C^{1,r}$ for some
$(1/2)<r<1$. As for (\ref{goal-2}), we note, as a preliminary step,
that
\begin{equation}\label{S-map}
{\mathcal S}\colon H^{-s}(\partial\Omega)\rightarrow
H^{-s+3/2}(\Omega)
\end{equation}
is well-defined and bounded for each $s\in [0,1]$, even when the
boundary of $\Omega$ is only Lipschitz. Indeed, with
$H^{-s+3/2}(\Omega)$ replaced by $H^{-s+3/2}(\Omega\cap B)$ for a
sufficiently large ball $B\subset\bbR^n$, this is proved in \cite{MT00} and
the behavior at infinity is easily taken care of by employing the
exponential decay of $E$.

For a fixed, arbitrary $j\in\{1,...,n\}$, consider next the operator
$\partial_{x_j}{\mathcal S}$ whose kernel is
$\partial_{x_j}E(x-y)=-\partial_{y_j}E(x-y)$. We write
\begin{equation}
\partial_{y_j}=\sum_{k=1}^n \nu_k(y)\nu_k(y)\partial_{y_j}
=\sum_{k=1}^n
\nu_k(y)\frac{\partial}{\partial\tau_{k,j}(y)}+\nu_j\partial_{\nu_y},
\end{equation}
where $\partial/\partial\tau_{k,j}=\nu_k\partial_j-\nu_j\partial_k$,
$j,k=1,\dots,n$, is a tangential derivative operator for which we have
\begin{equation}
\int_{\partial\Omega} d^{n-1}\sigma \, \frac{\partial
h_1}{\partial\tau_{j,k}}h_2
=-\int_{\partial\Omega}  d^{n-1}\sigma \, h_1\frac{\partial
h_2}{\partial\tau_{j,k}}, \quad h_1, h_2\in
H^{1/2}(\partial\Omega).
\end{equation}
It follows that
\begin{equation}\label{imp-id}
\partial_j{\mathcal S}h=-{\mathcal D}(\nu_j h)
+ \sum_{k=1}^n
{\mathcal S}\bigg(\frac{\partial(\nu_k h)}{\partial\tau_{k,j}}\bigg),
\end{equation}
where ${\mathcal D}$, the so-called double layer potential operator,
is the integral operator with integral kernel $\partial_{\nu_y}E(x-y)$.
Its mappings properties on the scale of Sobolev spaces have been analyzed
in \cite{MT00} and we note here that
\begin{equation}\label{D-map}
{\mathcal D}\colon H^{s}(\partial\Omega)\rightarrow
H^{s+1/2}(\Omega),\quad 0\leq s\leq 1,
\end{equation}
requires only that $\Omega$ is Lipschitz.

Assuming that multiplication by (the components of) $\nu$ preserves
the space $H^{1/2}(\partial\Omega)$ (which is the case if, e.g.,
$\Om$ is of class $C^{1,r}$ for some $(1/2)<r< 1$), the desired
conclusion about the operator (\ref{goal-2}) follows from (\ref{S-map}),
(\ref{imp-id}) and (\ref{D-map}). This concludes the proof of the fact that
$\dom(\wti H_{0,\Om}^N)\subseteq H^{2}(\Omega)$.

To prove that $\dom(\wti H_{0,\Om}^D)\subseteq H^{2}(\Omega)$ we proceed
in an analogous fashion, starting with the same representation
(\ref{sol}). This time, the requirement on $g$ is that
$Sg=h=\ga_D(\Delta-I_{\Omega})^{-1}f\in H^{3/2}(\partial\Omega)$, where
$S=\ga_D\circ{\mathcal S}$ is the trace of the single layer. Thus,
in this scenario, it suffices to know that
\begin{equation}\label{SS-iso}
S\colon H^{1/2}(\partial\Omega)\rightarrow H^{3/2}(\partial\Omega)
\end{equation}
is an isomorphism. When $\partial\Omega$ is of class $C^\infty$, it
has been proved in \cite[Proposition 7.9]{Ta96} that
$S\colon H^{s}(\partial\Omega)\to H^{s+1}(\partial\Omega)$ is an
isomorphism for each $s\in\bbR$ and, if $\Om$ is of class $C^{1,r}$ with
$(1/2)<r< 1$, the validity range of this result is limited to
$-1-r<s<r$, which covers (\ref{SS-iso}). The latter fact follows from an
inspection of Taylor's original proof of Proposition 7.9 in \cite{Ta96}.
Here we just note that the only significant difference is that if
$\partial\Omega$ is of class $C^{1,r}$ (instead of class $C^\infty$), then
$S$ is a pseudodifferential operator whose symbol exhibits a limited
amount of regularity in the space-variable. Such classes of operators have
been studied in, e.g., \cite{MMT01}, \cite[Chs.\ 1, 2]{Ta91}.
\end{proof}
%%%%%%%%%%%%%%%%%%%%%%%%%%%%%%%%%%%%%

We note that Lemma \ref{lA.1} also follows from \cite[Theorem 8.2]{DHP03}
in the case of $C^2$-domains $\Om$ with compact boundary. This is proved
in \cite{DHP03} by rather different methods and can be viewed as a
generalization of the classical result for bounded $C^2$-domains.

%%%%%%%%%%%%%%%%%%%%%%%%%%%%%%%%%%%%%
\begin{lemma}\label{lA.2}
Assume Hypothesis \ref{h6.1}\,$(i)$ and let %$q\in\bbR$. 
$q\in [0,1]$. Then for each
$z\in\bbC\backslash[0,\infty)$, one has
\begin{equation}\label{new6.45}
(H_{0,\Om}^D-zI_{\Om})^{-q},\,
(H_{0,\Om}^N-zI_{\Om})^{-q}\in\cB\big(L^2(\Om;d^nx),H^{2q}(\Om)\big).
\end{equation}
\end{lemma}
%%%%%%%%%%%%%%%%%%%%%%%%%%%%%%%%%%%%% 
\begin{proof}
For notational convenience, we denote by $H_{0,\Om}$ either one of the
operators $H_{0,\Om}^D$ or $H_{0,\Om}^N$. The operator $H_{0,\Om}$ is a
semibounded self-adjoint operator in $L^{2}(\Omega;d^nx)$, and thus the
resolvent set of $H_{0,\Om}$ is linearly connected.

Step $1$: We claim that it is enough to prove \eqref{new6.45} for one
point $z$ in the resolvent set of $H_{0,\Om}$. Indeed, suppose that
\eqref{new6.45} holds, and $z'$ is any other point in the resolvent
set of $H_{0,\Om}$. Connecting $z$ and $z'$ by a curve in the resolvent
set, and splitting this curve in small segments, without loss of
generality we may assume that $z'$ is arbitrarily close to $z$ so that
the operator $I_{\Omega}-(z'-z)(H_{0,\Om}-zI_{\Omega})^{-1}$ is invertible, and
thus the operator $(I_{\Omega}-(z'-z)(H_{0,\Om}-zI_{\Omega})^{-1})^{-q}$ is a
bounded operator on $L^2(\Om;d^nx)$. Then \eqref{new6.45} and the identity
\begin{equation}
(H_{0,\Om}-z'I_{\Omega})^{-q}=(H_{0,\Om}-zI_{\Omega})^{-q}
(I_{\Omega}-(z'-z)(H_{0,\Om}-zI_{\Omega})^{-1})^{-q}
\end{equation}
imply \eqref{new6.45} with $z$ replaced by $z'$, proving the claim.

Step $2$: By \cite[Theorem B.8]{Mc00} (cf.\ also Theorem
4.3.1.2 and Remark 4.3.1.2 in \cite{Tr78}), if
$\Omega\subseteq\bbR^n$ is a Lipschitz domain, $n\in\bbN$,
and $s_0,s_1\in\bbR$, then
\begin{equation}
\Big(H^{s_0}(\Omega),H^{s_1}(\Omega)\Big)_{\theta,2} =
H^{s}(\Omega), \quad  s=(1-\theta)s_0+\theta s_1,\;
0<\theta<1.
\end{equation}
Here, for Banach spaces $\mathcal{X}_0$ and $\mathcal{X}_1$, we denote by
$\big(\mathcal{X}_0,\mathcal{X}_1\big)_{\theta,p}$ the real
interpolation space (obtained by the $K$-method), as discussed, for
instance, in  \cite[Appendix B]{Mc00} and \cite[Sect.\ 1.3]{Tr78}. Letting
$s_0=0$, $s_1=2$, and $s=2q$, one then infers
\begin{equation}\label{INTERPOLATION}
\Big(L^{2}(\Omega;d^nx),H^{2}(\Omega)\Big)_{q,2}=H^{2q}(\Omega).
\end{equation}

Step $3$: Using the claim in Step $1$, we may assume without loss of
generality that $H_{0,\Om}-zI_{\Omega}$ is a strictly positive operator 
and thus the fractional power $(H_{0,\Om}-zI_{\Omega})^q$ can be defined
via its spectral decomposition (see, e.g., \cite[Sec.1.18.10]{Tr78}). We
remark that the operator
$(H_{0,\Om}-zI_{\Omega})^q$ is an isomorphism between the Banach space
$\dom(H_{0,\Om}-zI_{\Omega})^q$, equipped with the graph-norm, and the
space $L^{2}(\Omega;d^nx)$, and thus
\begin{equation}\label{BDD-q}
(H_{0,\Om}-zI_{\Omega})^{-q}\in
\cB\big(L^2(\Om;d^nx),\dom\big((H_{0,\Om}-zI_{\Omega})^q\big)\big).
\end{equation}
By an abstract interpolation result for strictly positive, self-adjoint
operators, see \cite[Theorem 1.18.10]{Tr78}, for any $\alpha,\beta\in\bbC$
with $\Re(\alpha),\Re(\beta)\ge0$ and $\theta\in(0,1)$ one has,
\begin{equation}\label{INTERPOLATION2}
\Big(\dom\big((H_{0,\Om}-zI_{\Omega})^\alpha\big),\dom\big((H_{0,\Om}
-zI_{\Omega})^\beta\big)
\Big)_{\theta,2}
=\dom\big((H_{0,\Om}-zI_{\Omega})^{\alpha(1-\theta)+\beta\theta}\big).
\end{equation}
Applying this result with $\alpha=0$ and $\beta=1$, one infers
\begin{equation}\label{INTERPOLATION3}
\Big(L^{2}(\Omega;d^nx),\dom(H_{0,\Om}-zI_{\Omega})\Big)_{q,2}
=\dom\big((H_{0,\Om}-zI_{\Omega})^q\big).
\end{equation}
Noting that $\dom(H_{0,\Om})=\dom(H_{0,\Om}-zI_{\Omega})$, and using
\eqref{INTERPOLATION},
\eqref{INTERPOLATION3}, and Lemma \ref{lA.1}, one arrives at the
continuous imbedding
\begin{equation}\label{INTERPOLATION4}
\dom\big((H_{0,\Om}-zI_{\Omega})^q\big) \hookrightarrow H^{2q}(\Omega).
\end{equation}
Thus, \eqref{new6.45} is a consequence of \eqref{BDD-q} and
\eqref{INTERPOLATION4}.
\end{proof}
%%%%%%%%%%%%%%%%%%%%%%%%%%%%%%%%%%%%%

Finally, we will prove an extension of a result of Nakamura \cite[Lemma
6]{Na01} from a cube in $\bbR^n$ to a Lipschitz domain $\Om$.
This requires some preparation. First, we note that
(\ref{A.16}) and (\ref{A.11a}) yield the following Green formula
\begin{align}
\langle\wti\gamma_N u,\ga_D\Phi\rangle = \big(\ol{\nabla u}, \nabla
\Phi\big)_{L^2(\Om;d^nx)^n} + \langle\Delta u,\Phi\rangle,
\lb{wGreen}
\end{align}
valid for any $u\in H^{1}(\Om)$ with $\Delta u\in\big(H^{1}(\Om)\big)^*$, and
any $\Phi\in H^{1}(\Om)$. The pairing on the left-hand side of (\ref{wGreen})
is between functionals in $\big(H^{1/2}(\dOm)\big)^*$ and elements in
$H^{1/2}(\dOm)$, whereas the last pairing on the right-hand side is
between functionals in $\big(H^{1}(\Om)\big)^*$ and elements in
$H^{1}(\Om)$. For further use, we also note that the adjoint of \eqref{6.1}
maps as follows
\begin{align}
\ga_D^* : \big(H^{s-1/2}(\dOm)\big)^* \to (H^{s}(\Om)\big)^*, \quad
1/2<s<3/2. \lb{ga*}
\end{align}

Next we observe that the operator $(\wti H^N_{0,\Om}-zI_\Om)^{-1}$,
$z\in\bbC\backslash\si(\wti H^N_{0,\Om})$, originally defined as
\begin{align}\label{fukcH}
(\wti H^N_{0,\Om}-zI_\Om)^{-1}:L^2(\Om;d^nx)\to L^2(\Om;d^nx),
\end{align}
can be extended to a bounded operator, mapping
$\big(H^{1}(\Om)\big)^*$ into $L^2(\Omega;d^nx)$. Specifically,
since $(\wti H^N_{0,\Om}-\ol{z}I_\Om)^{-1}: L^2(\Om;d^nx)
\to \dom(\wti H^N_{0,\Om})$ is bounded and since the inclusion
$\dom(\wti H^N_{0,\Om})\hookrightarrow H^{1}(\Om)$ is bounded, we can
naturally view $(\wti H^N_{0,\Om}-\ol{z}I_\Om)^{-1}$ as an operator
\begin{equation}
(\hatt H^N_{0,\Om}-\ol{z}I_\Om)^{-1} \colon L^2(\Om;d^nx) \to H^1(\Om)
\end{equation}
mapping in a linear, bounded fashion.
Consequently, for its adjoint, we have
\begin{align}\label{fukcH-bis}
\big((\hatt H^N_{0,\Om}-\ol{z}I_\Om)^{-1}\big)^* :
\big(H^{1}(\Om)\big)^* \to L^2(\Om;d^nx),
\end{align}
and it is easy to see that this latter operator extends the one
in \eqref{fukcH}. Hence, there is no ambiguity in retaining
the same symbol, that is, $(\wti H^N_{0,\Om}-zI_\Om)^{-1}$, both
for the operator in (\ref{fukcH-bis}) as well as for the operator
in (\ref{fukcH}). Similar considerations and conventions apply to
$(\wti H^D_{0,\Om}-zI_\Om)^{-1}$.

%%%%%%%%%%%%%%%%%%%%%%%%%%%%%%%%%%%%%
\begin{lemma} \lb{lA.3}
Let $\Om\subset\bbR^n$, $n\geq2$, be a Lipschitz domain and let
$z\in\bbC\backslash\big(\si(\wti H^D_{0,\Om})\cup\si(\wti H^N_{0,\Om})\big)$.
Then, on $L^2(\Omega;d^nx)$,
\begin{align} \lb{Na1}
(\wti H^D_{0,\Om}-zI_\Om)^{-1} - (\wti H^N_{0,\Om}-zI_\Om)^{-1} =
(\wti H^N_{0,\Om}-zI_\Om)^{-1}\ga_D^*\wti\gamma_N (\wti
H^D_{0,\Om}-zI_\Om)^{-1},
\end{align}
where $\ga_D^*$ is an adjoint operator to $\ga_D$ in the sense of
\eqref{ga*}
\end{lemma}
%%%%%%%%%%%%%%%%%%%%%%%%%%%%%%%%%%%%%
\begin{proof}
To set the stage, we note that the composition of operators appearing on
the right-hand side of \eqref{Na1} is meaningful since
\begin{align}
(\wti H^D_{0,\Om}-zI_\Om)^{-1} & \colon L^2(\Om;d^nx) \to
\dom(\wti H^D_{0,\Om}) \subset \{u\in H^1(\Om)\,|\,\Delta u\in (H^1(\Om))^*\},
\\
\wti\gamma_N & \colon \{u\in H^1(\Om)\,|\,\Delta u\in (H^1(\Om))^*\}
\to H^{-1/2}(\dOm)
\\
\ga_D^* & \colon \big(H^{1/2}(\dOm)\big)^*=H^{-1/2}(\dOm)
\to \big(H^{1}(\Om)\big)^*,
\\
(\wti H^N_{0,\Om}-zI_\Om)^{-1} & \colon \big(H^{1}(\Om)\big)^* \to
L^2(\Om;d^nx),
\end{align}
with the convention made just before the statement of the lemma
used in the last line.
Next, let $\phi_1,\psi_1\in L^2(\Om;d^nx)$ be arbitrary and define
\begin{align}
\begin{split}
\phi &= (\wti H^N_{0,\Om}-\ol{z}I_\Om)^{-1}\phi_1 \in \dom(\wti H^N_{0,\Om})
\subset H^{1}(\Om),
\\
\psi &= (\wti H^D_{0,\Om}-zI_\Om)^{-1}\psi_1 \in \dom(\wti H^D_{0,\Om})
\subset H^{1}(\Om).
\end{split} \lb{Na2}
\end{align}
It therefore suffices to show that the following identity holds:
\begin{align}
\begin{split}
&\big(\phi_1,(\wti H^D_{0,\Om}-zI_\Om)^{-1}\psi_1\big)_{L^2(\Om;d^nx)} -
\big(\phi_1,(\wti H^N_{0,\Om}-zI_\Om)^{-1}\psi_1\big)_{L^2(\Om;d^nx)}
\\
&\quad = \big(\phi_1,(\wti H^N_{0,\Om}-zI_\Om)^{-1}\ga_D^*\wti\gamma_N
(\wti H^D_{0,\Om}-zI_\Om)^{-1}\psi_1\big)_{L^2(\Om;d^nx)}.
\end{split}
\end{align}
We note that according to \eqref{Na2} one has,
\begin{align}
\big(\phi_1,(\wti H^D_{0,\Om}-zI_\Om)^{-1}\psi_1\big)_{L^2(\Om;d^nx)}
&= \big((\wti H^N_{0,\Om}-\ol{z}I_\Om)\phi,\psi\big)_{L^2(\Om;d^nx)},
\\
\big(\phi_1,(\wti H^N_{0,\Om}-zI_\Om)^{-1}\psi_1\big)_{L^2(\Om;d^nx)}
&= \big(\big((\wti H^N_{0,\Om}-zI_\Om)^{-1}\big)^*\phi_1,
\psi_1\big)_{L^2(\Om;d^nx)} \no
\\
&= \big((\wti H^N_{0,\Om}-\ol{z}I_\Om)^{-1}\phi_1,\psi_1\big)_{L^2(\Om;d^nx)}
\no
\\
&= \big(\phi,(\wti H^D_{0,\Om}-zI_\Om)\psi\big)_{L^2(\Om;d^nx)},
\end{align}
and, keeping in mind the convention adopted prior to the statement of the lemma,
\begin{align}
&\big(\phi_1,(\wti H^N_{0,\Om}-zI_\Om)^{-1}\ga_D^*\wti\gamma_N
(\wti H^D_{0,\Om}-zI_\Om)^{-1}\psi_1\big)_{L^2(\Om;d^nx)} \no
\\
&\quad =
\langle\ol{(\wti H^N_{0,\Om}-\ol{z}I_\Om)^{-1}\phi_1},\ga_D^*\wti\gamma_N
(\wti H^D_{0,\Om}-zI_\Om)^{-1}\psi_1\rangle \no
\\
&\quad = \big\langle\ol{\ga_D(\wti H^N_{0,\Om}-\ol{z}I_\Om)^{-1}\phi_1},
\wti\gamma_N (\wti H^D_{0,\Om}-zI_\Om)^{-1}\psi_1\big\rangle =
\big\langle\ol{\ga_D\phi},\wti\gamma_N\psi\big\rangle
\end{align}
where $\langle\cdot\,,\,\cdot\rangle$ stands for pairings
between Sobolev spaces (in $\Omega$ and $\partial\Omega$) and their duals.
Thus, matters have been reduced to proving that
\begin{align} \lb{Na3}
\big((\wti H^N_{0,\Om}-\ol{z}I_\Om)\phi,\psi\big)_{L^2(\Om;d^nx)} -
\big(\phi,(\wti H^D_{0,\Om}-zI_\Om)\psi\big)_{L^2(\Om;d^nx)} =
\big\langle\ol{\ga_D\phi},\wti\gamma_N\psi\big\rangle.
\end{align}
Using \eqref{wGreen} for the left-hand side of \eqref{Na3} one obtains
\begin{align}
&\big((\wti H^N_{0,\Om}-\ol{z}I_\Om)\phi,\psi\big)_{L^2(\Om;d^nx)} -
\big(\phi,(\wti H^D_{0,\Om}-zI_\Om)\psi\big)_{L^2(\Om;d^nx)} \no
\\
&\quad = -\big(\Delta\phi,\psi\big)_{L^2(\Om;d^nx)} +
\big(\phi,\Delta\psi\big)_{L^2(\Om;d^nx)}
\\
&\quad = \big(\nabla\phi,\nabla\psi\big)_{L^2(\Om;d^nx)^n} -
\big\langle\ol{\wti\gamma_N\phi},\ga_D\psi\big\rangle -
\big(\nabla\phi,\nabla\psi\big)_{L^2(\Om;d^nx)^n} +
\big\langle\ol{\ga_D\phi},\wti\gamma_N\psi\big\rangle \no
\\
&\quad = -\big\langle\ol{\wti\gamma_N\phi},\ga_D\psi\big\rangle +
\big\langle\ol{\ga_D\phi},\wti\gamma_N\psi\big\rangle. \no
\end{align}
Observing that $\wti\gamma_N\phi=0$ since $\phi\in\dom(H^N_{0,\Om})$,
one concludes \eqref{Na3}.
\end{proof}
%%%%%%%%%%%%%%%%%%%%%%%%%%%%%%%%%%%%%

%%%%%%%%%%%%%%%%%%%%%%%%%%%%%%%%%%%%%
\begin{remark} \lb{rA.4a}
While it is tempting to view $\ga_D$ as an unbounded but densely
defined operator on $L^2(\Om;d^nx)$ whose domain contains the space
$C_0^\infty(\Om)$, one should note that in this case its adjoint
$\ga_D^*$ is not densely defined: Indeed, the adjoint $\ga_D^*$ of $\ga_D$
would have to be an unbounded operator from $L^2(\dOm;d^{n-1}\si)$
to $L^2(\Om;d^nx)$ such that
\begin{equation}
(\ga_D f,g)_{L^2(\dOm;d^{n-1}\si)} = (f,\ga_D^* g)_{L^2(\Om;d^nx)}
\, \text{ for all }\, f\in\dom(\ga_D),\; g\in\dom(\ga_D^*).  \lb{A.60} 
\end{equation}
In particular, choosing $f\in C_0^\infty(\Om)$, in which case $\ga_D f =
0$, one concludes that $(f,\ga_D^* g)_{L^2(\Om;d^nx)}=0$ for
all $f\in C_0^\infty(\Om)$. Thus, one obtains $\ga_D^* g = 0$ for all
$g\in\dom(\ga_D^*)$. Since obviously $\gamma_D \neq 0$, \eqref{A.60}
implies $\dom(\gamma_D^*)=\{0\}$ and hence $\gamma_D$ is not a closable
linear operator in $L^2(\Om;d^nx)$. \hfill $\diamond$
\end{remark}
%%%%%%%%%%%%%%%%%%%%%%%%%%%%%%%%%%%%

%%%%%%%%%%%%%%%%%%%%%%%%%%%%%%%%%%%% 
\begin{remark} \lb{rA.4}
In the case of a domain $\Om$ of class $C^{1,r}$, $(1/2)<r<1$, the
operators $\wti H^D_{0,\Om}$ and $\wti H^N_{0,\Om}$
coincide with the operators $H^D_{0,\Om}$ and $H^N_{0,\Om}$,
respectively, and hence one can use the operators $H^D_{0,\Om}$ and
$H^N_{0,\Om}$ in Lemma \ref{lA.3}. Moreover, since $\dom(H^D_{0,\Om})
\subseteq H^2(\Om)$, one can also replace $\wti\gamma_N$ by $\ga_N$
(cf.\ \eqref{6.2}) in Lemma \ref{lA.3}. In particular,
\begin{align} \lb{Na1-bis}
(H^D_{0,\Om}-zI_\Om)^{-1} - (H^N_{0,\Om}-zI_\Om)^{-1}
=\big[\ga_D(H^N_{0,\Om}-\ol{z}I_\Om)^{-1}\big]^*\gamma_N
(H^D_{0,\Om}-zI_\Om)^{-1}.
\end{align}
${}$ \hfill $\diamond$
\end{remark}
%%%%%%%%%%%%%%%%%%%%%%%%%%%%%%%%%%%%%

%%%%%%%%%%%%%%%%%%%%%%%%%%%%%%%%%%%%%
\medskip
\noindent {\bf Acknowledgments.} We are indebted to Konstantin Makarov,
Alexander Pushnitski, Roland Schnaubelt, and Rico
Zacher for very helpful discussions.

Fritz Gesztesy and Yuri Latushkin gratefully acknowledge a
research leave for the academic year 2005/06 granted by the Research
Council and the Office of Research of the University of Missouri--Columbia.
Moreover, Yuri Latushkin gratefully acknowledges support by the Research
Board of the University of Missouri.
%%%%%%%%%%%%%%%%%%%%%%%%%%%%%%%%%%%%%%

%%%%%%%%%%%%%%%%%%%%%%%%%%%%%%%%%%%%%%
%%%%%%%%%%%%%%%%%%%%%%%%%%%%%%%%%%%%%% 

\end{document}